\documentclass[11pt]{article}

\def\mathbb{\Bbb}
\usepackage{amsfonts,latexsym,amstext}
\usepackage{amsmath}
\usepackage[active]{srcltx}
\usepackage{amssymb}

\newtheorem{theorem}{Theorem}[section]
\newtheorem{lemma}[theorem]{Lemma}
\newtheorem{proposition}[theorem]{Proposition}
\newtheorem{definition}{Definition}[section]

\newtheorem{hypothesis}[theorem]{Hypothesis}

\newtheorem{remark}[theorem]{Remark}
\newtheorem{corollary}[theorem]{Corollary}

\numberwithin{equation}{section}

\def\qed{{\hfill\hbox{\enspace${ \square}$}} \smallskip}
\def\sqr#1#2{{\vcenter{\vbox{\hrule height .#2pt \hbox{\vrule
 width .#2pt height#1pt \kern#1pt \vrule
width .#2pt} \hrule height .#2pt}}}}
\def\square{\mathchoice\sqr54\sqr54\sqr{4.1}3\sqr{3.5}3}

\def\ds{\begin{displaystyle}}
\def\eds{\end{displaystyle}}
\def\dis{\displaystyle }
\def\<{\langle }
\def\>{\rangle }

\def\dim{\noindent \hbox{{\bf Proof.} }}

\def\R{\mathbb R}

\def\E{\mathbb E}
\def\P{\mathbb P}
\def\Q{\mathbb Q}

\def\calf{{\cal F}}
\def\calg{{\cal G}}

\def\calm{{\cal M}}
\def\calp{{\cal P}}

\def\call{{\cal L}}
\def\cals{{\cal S}}
\def\calo{{\cal O}}

\begin{document}

\title{A Bismut-Elworthy formula for quadratic BSDEs}
\date{}
 \author{
 Federica Masiero\\
 Dipartimento di Matematica e Applicazioni, Universit\`a di Milano Bicocca\\
 via Cozzi 55, 20125 Milano, Italy\\
 e-mail: federica.masiero@unimib.it
}

\maketitle  
\begin{abstract}
We consider a backward stochastic differential equation in a Markovian framework for the pair of processes $(Y,Z)$, with generator with quadratic growth with respect to $Z$. Under non-degeneracy assumptions, we prove an analogue of the well-known Bismut-Elworty formula when the generator has quadratic growth with respect to $Z$. Applications to the solution of a semilinear Kolmogorov equation for the unknown $v$ with nonlinear term with quadratic growth with respect to $\nabla v$ and final condition only bounded and continuous are given, as well as applications to stochastic optimal control problems with quadratic growth.
\end{abstract}

\section{Introduction}
In this paper we study a Bismut-Elworthy type formula for BSDEs in a Markovian framework
when the generator has quadratic growth with respect to $Z$.

\noindent Namely, let us consider a forward stochastic differential equation in the Hilbert space $H$,
\begin{equation}
\left\{
\begin{array}
[c]{l}%
dX^{t,x}_\tau  =AX^{t,x}_\tau d\tau+F(\tau, X^{t,x}_\tau)+G(\tau)dW_\tau
,\text{ \ \ \ }\tau\in\left[  t,T\right] \\
X^{t,x}_t =x,
\end{array}
\right.  \label{forwardintro}%
\end{equation}
where $A$ is the generator of a strongly continuous semigroup in $H$, $\left\lbrace W_\tau,\,\tau\geq 0\right\rbrace$
is a cylindrical Wiener process in $H$, and the maps 
 $F$ and $G$ are defined on $:[0,T]\times H$ and on $[0,T]$ respectively, and take their values in $H$ and $L(H;H)$, respectively.
The solution of equation (\ref{forwardintro}) will be denoted by $X_\tau$, or also by $X_\tau^{t,x}$ to stress
the dependence on the initial conditions. The transition semigroup related to $X^{t,x}$ is denoted by
\[
 P_{t,\tau}[\phi](x):=\E\phi(X_\tau^{t,x}).
\]
At least formally, the generator of $P_{t,\tau}$ is the second order differential operator
\[
 (\call_t f)(x)=\frac{1}{2}(Tr G(t)G^*(t) \nabla^2 f)(x)+\<Ax,\nabla f(x)\>+\<F(t,x),\nabla f(x),
\quad0\leq t\leq T\>.
\]
Under the invertibility assumption on $G$,
\[
\vert G^{-1}(\tau)\vert \leq C,
\]
and in the case of $G$ depending also
on $x\in H$, the infinite dimensional extension for the Bismut-Elworthy formula has been proved, see e.g. \cite{Ce}, \cite{DP3}.
According to the Bismut formula, for every $0\leq t<\tau\leq T,\,x\in H$, for every direction $h\in H$,
and for every bounded and continuous real function $f$ defined on $H$,
 \begin{equation}\label{Bismutintro}
 \<\nabla_x P_{t,\tau}[f](x), h\>=\E f\left(X_\tau^{t,x}\right)U^{h,t,x}_\tau,
\end{equation}
where
\begin{equation*}
 U^{h,t,x}_\tau:=\dfrac{1}{\tau-t}\int_t^\tau \<G^{-1}(r,X^{t,x}_r)\nabla_xX_r^{t,x}h,dW_r \>
\end{equation*}
Also in \cite{Ce} and \cite{DP3},
the Bismut formula has been used as a basic tool to
study the semilinear Kolmogorov equation driven by $\call_t$:
\begin{equation}
\left\{
\begin{array}
[c]{l}%
-\frac{\partial v}{\partial t}(t,x)=\call_t v\left(  t,x\right)
+\psi\left( t,x,v(t,x),\nabla v(t,x)G(t)   \right)  ,\text{ \ \ \ \ }t\in\left[  0,T\right]
,\text{ }x\in H\\
v(T,x)=\phi\left(  x\right)  ,
\end{array}
\right.  \label{Kolmointro}
\end{equation}
in the case of hamiltonian function $\psi$ lipschitz continuous with respect to $\nabla v$. We point out that the hamiltonian function $\psi$, here in equation (\ref{Kolmointro}), and throughout the paper, depends  on the directional derivative $\nabla v G$, but since
$G$ is invertible this is equivalent to consider an hamiltonian function depending on the derivative $\nabla v$.

 We recall that by mild solution for the Kolmogorov equation (\ref{Kolmointro}) we mean a bounded and continuous function $v:[0,T]\times H \rightarrow H$, once G\^ateaux differentiable with respect to $x$,
such that $v$ satisfies the integral equality
\begin{equation} 
v(t,x)=P_{t,T}\left[  \phi\right]  \left(  x\right)  +\int_{t}^{T}%
P_{t,s}\left[  \psi(s,\cdot, v(s,\cdot), \nabla v\left(  s,\cdot\right)G(s)) 
\right]  \left(  x\right)  ds.
\text{\ \ }t\in\left[  0,T\right]  ,\text{ }x\in H. \label{solmildkolmointro}%
\end{equation}
Notice that as a byproduct of the Bismut formula (\ref{Bismutintro}), we get the estimate
\[
\vert\<\nabla_x P_{t,\tau}[f](x), h\>\vert \leq C (\tau-t)^{-\frac{1}{2}}\vert h\vert\Vert f\Vert_\infty.
\]
Due to the lipschitz character of $\psi$, this is sufficient to prove, by a fixed point argument, existence and
uniqueness of a mild solution of the Kolmogorov equation (\ref{Kolmointro}), see e.g. \cite{Go1}.

In the present paper we aim to solve equation (\ref{Kolmointro}) when $\psi$ has quadratic growth with respect
to $\nabla v$: in this case a Bismut formula for the directional derivative of the transition semigroup is no more sufficient to
solve the Kolmogorov equation (\ref{Kolmointro}) by a fixed point argument, see \cite{Go2} and \cite{Mas1} where the
solution of the Kolmogorov equation in the quadratic case is performed.

\noindent To this aim, we want to prove a Bismut-Elworthy formula in a nonlinear situation, following \cite{futeBismut}.
Namely, we recall that connections between partial differential equations and stocahstic differential equations of backward type ( BSDEs in the following ) have been established in the pioneering papers \cite{PaPe1}, \cite{PaPe},
 in the finite domensional case, and have been extended to the infinite dimensional case in \cite{fute}. More precisely, let us consider the forward backward system
\begin{equation}
 \label{fbsdeintro}
    \left\{\begin{array}{l}\dis dX^{t,x}_\tau =
AX^{t,x}_\tau d\tau+ F(\tau, X^{t,x}_\tau)+G(\tau)dW_\tau,\quad \tau\in
[t,T]\subset [0,T],
\\\dis
X^{t,x}_t=x,
\\\dis
 dY_\tau^{t,x}=-\psi(\tau,X^{t,x}_\tau,Y^{t,x}_\tau,Z^{t,x}_\tau)\;d\tau+Z^{t,x}_\tau\;dW_\tau,
  \\\dis
  Y_T^{t,x}=\phi(X_T^{t,x}),
\end{array}\right.
\end{equation}
where the forward equation is just equation (\ref{forwardintro}), and the generator $\psi$ and the final condition
$\phi$ in the BSDE are just the nonlinear term ( also referred to as hamiltonian ) and the final condition in the semilinear Kolmogorov equation (\ref{Kolmointro}).
If we set $(X^{t,x}, Y^{t,x},Z^{t,x})$ the solution of the FBSDE (\ref{fbsdeintro}), it is well known that
the function 
\[
 v(t,x):=Y_t^{t,x}
\]
is a mild solution of the Kolmogorov equation (\ref{Kolmointro}), and moreover
the following identification for $Z$ holds true
\[
 \nabla v(t,x)= Z^{t,x}_tG^{-1}(t), \quad t\in [0,T].
\]
In the present paper we aim to solve a semilinear Kolmogorov equation like (\ref{Kolmointro})
with $\psi$ quadratic with respect to $\nabla v(t,x) G(t)$, and to do this we want to extend to our context the
following nonlinear version of the Bismut-Elworthy formula, see \cite{futeBismut}: for $0\leq t<\tau\leq T,\,x\in H$, for every direction $h\in H$,
\begin{equation}\label{Bismut-nonlin-intro}
 \E\left[ \nabla_x\,Y^{t,x}_sh \right]=
\E\int_s^T\psi\left(r,X_r^{t,x},Y_r^{t,x},Z_r^{t,x}\right)U^{h,t,x}_r\,dr
+\E\left[  \phi(X_T^{t,x})U^{h,x}_T\right]
\end{equation}
for $\phi $ and $\psi$ bounded and continuous functions satisfiying suitable assumptions, including lipschitzianity of $\psi$ with respect to $Z$.
Here we prove a nonlinear Bismut formula like (\ref{Bismut-nonlin-intro}) when $\psi$ has quadratic growth with respect to $Z$.
One of the main tools in proving this formula is the apriori estimates on $Z^{t,x}$, where 
$(X^{t,x},Y^{t,x},Z^{t,x})$ is solution to the FBSDE (\ref{fbsdeintro}). Namely,
with techniques similar to the ones used in \cite{BaoDeHu} and \cite{Ri},
we are able to prove that
\begin{equation}\label{stimabismutintro}
\vert Z^{t,x}_t\vert\leq C (T-t)^{-1/2},
\end{equation}
where $C$ depends on $t,\;T,\;A,\;F,\;\Vert\phi\Vert_\infty$. We remark that we are able to prove the
fundamental apriori estimates (\ref{stimabismutintro}) only in the case of $G$ not depending on $x$, and consequently the Bismut formula in the quadratic case is proved under this restrictive assumption which is not present in the lipschitz continuous case considered in \cite{futeBismut}, see also remark \ref{remark:prop-aprioriZ} for further and more technical comments on this point.

\noindent If the coefficients are differentiable, beside estimates (\ref{stimabismutintro}),
 we can also prove, as in \cite{Ri1} and \cite{MR}, that
\begin{equation*}
\vert Z^{t,x}_t\vert\leq C ,
\end{equation*}
where $C$ depends on the derivatives of $\phi$ and $\psi$.
This allows to prove a nonlinear Bismut formula in the quadratic case for differentiable $\psi$ and $\phi$,
and by an approximation procedure differentiability assumptions can be removed.

The nonlinear Bismut formula (\ref{Bismut-nonlin-intro}), which has its own idependent interest, allows to solve the Kolmogorov equation with hamiltonian function quadratic with respect to $\nabla v$. As we have previously discussed, the solution of the Kolmogorov equation (\ref{Kolmointro}) in the quadratic case is not a consequence of the linear Bismut formula.

Second order differential equations are a widely studied topic in the literature, see e.g. \cite{DP3}.
 In particular mild solutions of semilinear Kolmogororv equations with the structure of equation
(\ref{Kolmointro}) and with $\psi$ lipschitz continuous are studied both by an
analytic approach, see e.g. \cite{Go1} and \cite{Mas}, by a purely probabilistic approach,
by means of backward stochastic
differential equations (BSDEs in the following), see \cite{fute}.

In the case of $\psi$ only locally lipschitz continous, we cite \cite{Go2}, \cite{MR} and \cite{Mas1}, where in
particular the quadratic case is studied with datum $\phi$ only continuous.
In the present paper we consider an hamiltonian function locally lipschitz continuous
and with quadratic growth with respect to $z$, and with respect to $x$ we ask the hamiltonian function $\psi$
and the final datum $\phi$ to be bounded, and no further regularity than continuity in $x$ is asked.
The case of an hamiltonian function which is quadratic and locally lipschitz continuous with respect to $z$ is addressed also in \cite{Mas1} with the same assumptions of final datum $\phi$ bounded and continuous with respect to $x$, but with different asumptions on $G$: here $G$ is assumed to be invertible, while in \cite{Mas1}
$G$ and $A$ commute.
In \cite{Mas1} it is also addressed the case of locally lipschitz continuous hamiltonian function with superquadratic growth with respect to $z$. In this case the final datum is assumed to be bounded and lipschitz continuous with respect to $x$. The paper \cite{MR} is a generalization in this direction: the hamiltonian function is locally lipschitz continuous with respect to $z$, and it is allowed also superquadratic growth with respect to $z$, with respect to $x$ the hamiltonian function and the final datum are allowed to have polynomial growth, but they are assumed to be locally lipschitz continuous
with respect to $x$. See also remark \ref{remarkKolmo} and \ref{remarkKolmoheat} for further technical comments.

\noindent We also cite the paper \cite{BriFu} where quadratic infinite dimensional
HJB equations are solved
by means of BSDEs: the generator $\mathcal{L}$ is related to
a more general Markov process $X$ than the one considered here in
(\ref{forwardintro}), and no assumptions on the diffusion
coefficient are made,
but only the case of final condition $\phi$ and generator $\psi$ G\^ateaux differentiable is treated.

In the present paper, by applying the non linear Bismut formula (\ref{Bismut-nonlin-intro}), we are able to prove
existence and uniqueness of a mild solution for the Kolmogorov equation (\ref{Kolmointro})
with quadratic hamiltonian function, lipschitz continuous with respect to $x$, 
and locally lipschitz continuous with respect to $\nabla v$, and with final condition
only bounded and continuous.

The results are applied to a stochastic optimal control problem, related to a controlled state equation
\begin{equation}
\left\{
\begin{array}
[c]{l}%
dX^{u}_\tau  =AX^{u}_\tau d\tau +F(\tau,X_\tau)d\tau+ u_\tau  d\tau+
G(\tau)dW_\tau ,\text{ \ \ \ }\tau\in\left[  t,T\right] \\
X^{u}_t  =x.
\end{array}
\right.  \label{sdecontrolintro}%
\end{equation}
and to a cost functional
\begin{equation*}
J\left(  t,x,u\right)  =\mathbb{E}\int_{t}^{T}
l\left(s,X^{u}_s,u_s\right)ds+\mathbb{E}\phi\left(X^{u}_T\right). 
\end{equation*}
where $l$ has quadratic growth with respect to $u$, and the admissible controls $u$
are not asked to take values in a bounded set.

The paper is organized as follows:
in section \ref{sezionesde} some results on the forward equation are collected,
in section \ref{sez-fbsde} connections between BSDEs and Kolmogorov equations
are recalled, and the Bismut-Elworthy formula proved in \cite{futeBismut} is presented.
Section \ref{sez-Bismut-quad} deals with the Bismut-Elworthy formula in the quadratic case, and
with $\psi$ smooth and $\phi$ only continuous, in section \ref{sez-BismutPDE}, starting from the Bismut formula proved in \ref{sez-Bismut-quad}, the Kolmogorov equation
(\ref{Kolmointro}) is solved with $\psi$ quadratic and only lipschitz continuous
and $\phi$ only continuous, and in \ref{sez-appl-contr} applications to control are given, in particular in
 \ref{sez_contr_heat} the results are applied to a
controlled heat equation.

\section{Notations and preliminary results on the forward equation}
\label{sezionesde}

Throughout the paper, we let $(\Omega, \calf, \P)$ be a complete probability
space and we denote by $H$ and $\Xi$ real and separable Hilbert spaces. We consider a cylindrical
Wiener process $(W_t)_{t\geq 0}$ with values in $\Xi$, and defined on $(\Omega, \calf, \P)$. For $t\geq 0$, let $\calf_t$ denote the $\sigma$-algebra
generated by $(W_s,\, s\leq t)$ and augmented with the $\P$-null sets of $\calf$. The notation
$\E_t$ stands for the conditional expectation given $\calf_t$.

\noindent For any real and separable Hilbert space $K$, we denote further
\begin{itemize}
 \item $\cals^p(K), \, 1\leq p<\infty$, or $\cals^p$ where no confusion is possible,
the space of all predictable processes
$(Y_t)_{t\in[0,T]}$ with values in $K$, normed by 
$$\Vert Y\Vert _{\cals^p}=\left(\E\sup_{t\in[0,T]}\vert Y_t\vert^p\right)^{1/p};$$
 $\cals^\infty(K)$,
or $\cals^\infty$ where no confusion is possible, the space of all bounded predictable processes.
\item $\calm^p(K),\, 1\leq p<\infty$, or $\calm^p$ where no confusion is possible, the space of all predictable processes
$(Z_t)_{t\in[0,T]}$ with values in $K$, normed by 
$$\Vert Z\Vert _{\calm^p}=\left(\E\left(\int_0^T\vert Z_t\vert^2dt\right)^{p/2}\right)^{1/p}.$$
\end{itemize}

Following \cite{fute}, given two Banach spaces $E$ and $V$ we say that a function
$f:X\rightarrow V $ belongs to the class $\mathcal{G}^{1}\left(  E, V\right)  $ if$\ f$ is continuous
and G\^ateaux differentiable on $E$, and the gradient $\nabla f:E\rightarrow
L\left(  E,V\right)  $ is strongly continuous, that is for every
directions $e\in E$ the map $\nabla f\left(  \cdot\right)  e:E\rightarrow V$ is continuous.
Generalizations of this definition for functions depending on several variables are immediate, see also \cite{fute}.

\subsection{The forward equation}
\label{subsez-sde}
We consider the Markov process $X$ (also denoted $X^{t,x}$ to stress the dependence
on the initial conditions) in $H$ solution to equation%
\begin{equation}
\left\{
\begin{array}
[c]{l}%
dX^{t,x}_\tau  =AX^{t,x}_\tau d\tau+F(\tau, X^{t,x}_\tau)+G(\tau)dW_\tau
,\text{ \ \ \ }\tau\in\left[  t,T\right] \\
X^{t,x}_t =x,
\end{array}
\right.  \label{forward}%
\end{equation}
where $(W_\tau)_{\tau\in [0,T]}$ is a cylindrical Wiener process with values in $\Xi$.
On the coefficients of equation (\ref{forward}) we assume the following:
\begin{hypothesis}
 \label{ip_forward}
\begin{enumerate}
 \item The linear operator $A$ is the generator of a strongly continuous
semigroup $\left(  e^{t A},t\geq0\right)  $ in the Hilbert space $H.$
\item The map $F:[0,T]\times H\rightarrow H$ is measurable and satisfies, for some constant $C>0$, 
\begin{align*}
& \vert F(t,x)  \vert\leq C (1+\vert x\vert),\\ \nonumber
& \vert F(t,x)  -F(t,y)\vert\leq C \vert x-y\vert,
\end{align*}
for every $t\in[0,T],\,x, y\in H$.
\item The map $G:[0,T]\times H\rightarrow L(\Xi,H)$ is such that $\forall\, \xi\in\Xi$
the map $G\xi:[0,T] H\rightarrow H$ is measurable; for every
$s>0,\,t\in[0,T]$ $e^{sA}G(t,x)\in L_2(\Xi,H)$ and
the following estimate hold true
 \begin{align*}
 & \vert e^{sA}G(t)  \vert_{ L_2(\Xi,H)}\leq C s^{-\gamma},\\ \nonumber
  \end{align*}
for some constant $C>0$ and $0\leq\gamma< \dfrac{1}{2}.$
\end{enumerate}
\end{hypothesis}
We need further to assume, on the coefficients $F$ and $G$, the following
\begin{hypothesis}
 \label{ip_forward_agg}
\begin{enumerate}
\item For every $s>0$ and $t\in[0,T]$,
\begin{equation*}
 F(t,\cdot) \in \calg^1(H,H);
\end{equation*}
\item for some constant $C>0$ 
\begin{equation*}
\vert G(t)\vert_{L(\Xi,H)}\leq C,\qquad t\in[0,T].
\end{equation*}
\end{enumerate}
\end{hypothesis}
For our main result, we will further assume later that the operators $G(t)$ are boundedly invertible.
\begin{proposition}
 \label{prop-ex-forward}
Under hypothesis \ref{ip_forward}, for every $p\geq 2$, there exists a unique process $X\in \cals^p$
solving equation (\ref{forward}) and satisfying moreover
\[
 \E\sup_{\tau\in[0,T]}\vert X^{t,x}_\tau\vert^p\leq C(1+\vert x\vert^p).
\]

\noindent If also hypothesis \ref{ip_forward_agg} holds true, we get that the map
$(t,x)\mapsto X^{t,x}$ belongs to $\calg^{0,1}([0,T]\times H,\cals^p), \,\forall 1\leq p\leq \infty$ 
and for every direction $h\in H$, the directional derivative $\nabla_x X_\tau^{t,x}h$ solves
the following equation
\begin{equation}
\nabla_x X_\tau^{t,x}h =e^{(\tau-t)A}h+\displaystyle\int_t^\tau e^{(\tau-\sigma)A}\nabla_xF(\sigma, X_\sigma^{t,x})\nabla_x X_\sigma^{t,x}h\,d\sigma
 \label{forward-der}%
\end{equation}
\end{proposition}
\dim The proof follows by \cite{DP1} for the existence part, and by \cite{fute}
for the differentiability part. Notice that since $G$ does not depend on $x$, by a simple application of
the Gronwall lemma we get
\begin{equation}\label{stima-nablaX}
\vert \nabla_x X_\sigma^{t,x}h\vert \leq C \vert h\vert 
\end{equation}
where $C$ is uniform with respect to $\sigma $ and may depend on $t,\, T,\, A$ and
on the lipschitz constant of $F$.
\qed
\section{The forward-backward system and connections with PDEs}
\label{sez-fbsde}
In this section we consider the following 
forward-backward system: for given $t\in [0,T]$ and $x\in H$,
\begin{equation}\label{fbsde}
    \left\{\begin{array}{l}\dis dX_\tau =
AX_\tau d\tau+ F(\tau, X_\tau)+G(\tau)dW_\tau,\quad \tau\in
[t,T]\subset [0,T],
\\\dis
X_t=x,
\\\dis
 dY_\tau=-\psi(\tau,X_\tau,Y_\tau,Z_\tau)\;d\tau+Z_\tau\;dW_\tau,
  \\\dis
  Y_T=\phi(X_T),
\end{array}\right.
\end{equation}
for the unknown $(X,Y,Z)$, also denoted by
$(X^{t,x},Y^{t,x},Z^{t,x})$ to stress the dependence on the initial
conditions $t$ and $x$. The process $X$, which is solution of the forward equation (\ref{forward}),
has been extended for $0\leq s\leq t$ by setting
$X_s=x$ for $0\leq s\leq t$. 
The second equation is of backward type for the unknown $(Y,Z)$
and depends on the Markov process $X$.
Under suitable assumptions on the coefficients
 $\psi:[0,T]\times H\times \R\times\Xi
\rightarrow\mathbb{R}$
and $\mathbb{\phi}:H\rightarrow\mathbb{R}$
we will look for a solution consisting of a pair of predictable processes,
taking values in $\mathbb{R}\times \Xi$, such that $Y$ has
continuous paths and
\[
\|\left( Y,Z\right)\|^2_{\cals^2\times\calm^2}:=
\mathbb{E}\sup_{\tau\in\left[ 0,T\right] }\left\vert Y_{\tau}\right\vert
^{2}+\mathbb{E}\int_{0}^{T}\left\vert Z_{\tau}\right\vert ^{2}d\tau<\infty,
\]
see e.g. \cite{PaPe1} for the classical starting case where the generator $\psi$
is assumed to be lipschitz continuous with respect to $Y$ and $Z$.
In the present paper we assume that the generator is lipschitz continuous
with respect to $y$ and locally lipschitz continuous with respect to $z$,
namely we assume that with respect to $z$ the generator $\psi$ has quadratic growth,
as stated in the following:

\begin{hypothesis}
\label{ip-psiphi} The function $\phi$ is continuous and the function $\psi$ is measurable, moreover for every fixed $t\in[0,T]$ the map $\psi(t,\cdot,\cdot,\cdot):H\times\R\times\Xi\rightarrow\R$
is continuous. There exist nonnegative constants $L_\psi,\,K_\psi,\,K_\phi$
such that
\begin{align*}
&\vert \psi(t,x_1,y_1,z_1)- \psi(t,x_2,y_2,z_2)\vert\leq L_\psi\left(\vert x_1-x_2\vert+\vert y_1-y_2\vert+\vert z_1-z_2\vert
(1+\vert z_1\vert+\vert z_2\vert)\right),\\
&\vert \psi(t,x,0,0)\vert\leq K_\psi, \qquad 
\vert \phi(x)\vert\leq K_\phi,
\end{align*}
for every $t\in[0,T]$, $x_1,x_2\in H$, $y_1,y_2\in\R$ and $z_1,z_2\in \Xi$
\end{hypothesis}

\begin{theorem}\label{teo-ex-bsdequadr}Let $(X,Y,Z)$ be solution of the forward-backward system (\ref{fbsde}),
and assume that hypotheses \ref{ip_forward} and \ref{ip-psiphi} hold
true.
Then there exists a unique solution of the markovian BSDE in (\ref{fbsde}) such that 
\[
 \Vert Y\Vert _{\cals^2}+\Vert Z\Vert_{\calm^2}\leq C,
\]
where $C$ is a constant that may depend on $A,\,F,\,G,\, K_\psi,\,L_\psi,\,K_\phi, t,\,T$ but not on $x$.
\end{theorem}
\dim This result substantially follows from \cite{Kob}.

\qed

Many other results concerning not bounded final data and generator with polynomial growth with respect to $x$
have been proved after \cite{Kob}, we cite here \cite{BriHu2006} and \cite{BriHu2008}.
In the Markovian framework, if $G(t)=G$ and if $\mu<2$, and moreover if for every $x,x'\in H$
$\phi$ and $\psi$ satisfy
\[
 \vert \psi(t,x,y,z)-\psi(t,x',y,z)\vert \leq 
\left(C+\frac{\beta}{2}\vert x\vert^r+\frac{\beta}{2}\vert x'\vert^r\right)\vert x -x'\vert;
\]
 and 
\[
 \vert \phi(x)-\psi(x')\vert \leq 
\left(C+\frac{\alpha}{2}\vert x\vert^r+\frac{\alpha}{2}\vert x'\vert^r\right)\vert x -x'\vert;
\]
then it has been proved in \cite{Ri1} existence and uniqueness of a solution for a markovian BSDE,
like the one considered in \ref{fbsde}. Namely in \cite{Ri} it has been proved that there exists a unique solution of the markovian BSDE in (\ref{fbsde}) such that 
\[
 \Vert Y\Vert _{\cals^2}+\Vert Z\Vert_{\calm^2}<C(1+\vert x\vert^{\mu}) 
\text{ and }\vert Z^{t,x}\vert\leq C\left( 1+\vert X^{t,x}\vert ^r\right).
\]
In \cite{MR}, the extension to the case of $X$
taking values in an infinite dimensional Hilbert space is considered 

In the present paper we consider the case of $\phi$ bounded, continuous
without lipshitz properties, and $\psi$ bounded and lipschitz continuous
with respect to $x$, and with quadratic growth with respect to $z$.

We also cite the Feynman-Kac formula, proved in \cite{BriFu} when all the coefficients are differentiable
and in the case of $\psi$ quadratic with respect to $z$, and generalized e.g. in \cite{Mas1} to nonsmooth coefficients, and in \cite{MR} to the case of $\psi$ superquadratic with respect to $z$.
More precisely, let $\mathcal L $ be the generator of the transition
semigroup $(P_{t,s})_{0\leq t\leq s\leq T}$, that is, at least formally,
$$
(\call_t f)(x)=\frac{1}{2}(Tr G(t)G^*(t) \nabla^2 f)(x)+\<Ax,\nabla f(x)\>+\<F(t,x),\nabla f(x)\>.
$$
Let us consider the following equation
\begin{equation}
\left\{
\begin{array}
[c]{l}%
\frac{\partial v}{\partial t}(t,x)=-\call_t v\left(  t,x\right)
+\psi\left( t,x,v(t,x),\nabla v(t,x)G(t)   \right)  ,\text{ \ \ \ \ }t\in\left[  0,T\right]
,\text{ }x\in H\\
v(T,x)=\phi\left(  x\right)  ,
\end{array}
\right.  \label{Kolmo}%
\end{equation}
We introduce the notion of mild solution of the non linear Kolmogorov
equation (\ref{Kolmo}), see e.g. \cite{fute}. 
Let $P_{t,\tau},\,t\leq\tau\leq T$,
the transition semigroup related to the process $X^{t,x}$ solution of the forward equation (\ref{forward}),
namely, for every bounded and measurable function $\phi:H\rightarrow\R$
\[
 P_{t,\tau}[\phi](x)=\E\phi(X_\tau^{t,x}).
\]
Since $\call_t$ is at least formally the generator of
$(P_{t,s})_{0\leq t\leq s\leq T}$, the variation of constants formula for (\ref{Kolmo}) gives:%
\begin{equation}
v(t,x)=P_{t,T}\left[  \phi\right]  \left(  x\right)  +\int_{t}^{T}%
P_{t,s}\left[  \psi(s,\cdot, v(s,\cdot), \nabla v\left(  s,\cdot\right)G(s,\cdot)) 
\right]  \left(  x\right)  ds.
\text{\ \ }t\in\left[  0,T\right]  ,\text{ }x\in H. \label{solmildkolmo}%
\end{equation}
We use this formula to give the notion of mild
solution for the non linear Kolmogorov equation (\ref{Kolmo}).

\begin{definition}
\label{defsolmildkolmo}A function $v:\left[  0,T\right]  \times H\rightarrow\mathbb{R}$ is a mild
solution of the non linear Kolmogorov equation (\ref{Kolmo}) if $v\in C_{b}\left(  \left[  0,T\right]  \times H\right)  $ and it is differentiable, namely $v\in\calg^{0,1}\left(  \left[  0,T\right]  \times H\right)$
and equality (\ref{solmildkolmo}) holds.
\end{definition}
We are ready to give a precise statement of the Feynman Kac formula in the quadratic case.
\begin{theorem}\label{teo_fey_kac}Let hypotheses \ref{ip_forward}, \ref{ip_forward_agg} and \ref{ip-psiphi}
hold true.  Moreover assume that $\phi$ is G\^ateaux differentiable with a bounded derivative, and that $\psi$
is G\^ateaux differentiable with respect to $x$, $y$ and $z$ with bounded derivatives. The nonlinear Kolmogorov equation (\ref{Kolmo}) has a unique mild solution $v$ given by the formula
$$v(t, x) = Y_t^{t,x} ,\qquad     (t, x) \in [0, T ] × H$$
where $(X^{t,x},Y ^{t,x} , Z^{ t,x} )$ is the solution to the FBSDE \ref{fbsde}. Moreover, we
have, $\P$-a.s.,
\[
Y_s^{t,x} = v(s, X_s^{t,x} ),\quad       Z_s^{t,x} =\nabla_x v(s, X_s^{t,x} )\nabla_x  X_s^{t,x}G(s).
\]
Moreover under our assumptions, there exists a constant $C$, that may depend also on $\nabla_x\phi$,
$\nabla_x\psi$ and $L_\psi$, such that
\begin{equation}\label{stimaZ-diffle}
\vert Z_s^{t,x}\vert \leq C
\end{equation}
\end{theorem}
\dim
The result is proved in \cite{BriFu}, Proposition 12 and Theorem 15, also for the case of $G$ depending on $x$.
\qed

We now prove an apriori estimate on $Z^{t,x}$ depending only on the growth of final datum $\phi$
with respect to $x$; since in our framework the coefficients are bounded with respect to $x$,
this apriori estimates actually depends on the $L^{\infty}$-norm of the final datum. To this aim,
we need to formulate our fundamental assumption on the invertibility of $G$.
\begin{hypothesis}\label{ip_G_invertible} For every $t\in[0,T]$, the operator $G(t)$
has a bounded inverse and there exists a constant $B$ such that
 \[
\vert G^{-1}(t)\vert \leq B,\qquad t\in[0,T].
 \]
\end{hypothesis}
\begin{remark}\label{remark:prop-aprioriZ}
In the following proposition we prove the fundamental apriori estimate on $Z$ (\ref{stimabismut}): we prove it by BSDEs' techniques and we need here that $G$ does not depend on $x$. We are aware that in finite dimension  a similar estimates in proved in \cite{CriDel} for operators $G$ depending also on $x$. The extension of this estimate to infinite dimensions is an open problem up to our knowledge. We also notice that in the case of lipwschitz generator $\psi$ treated in \cite{futeBismut} such an apriori estimate on $Z$ is not necessary to prove the Bismut-Elworthy formula, while it is a fundamental tool in our proof of the Bismut formula in the quadratic case.
\end{remark}
\begin{proposition}\label{prop-aprioriZ}
 Let $(Y,Z)$ be the solution of the BSDE in (\ref{fbsde}). Let 
hypotheses \ref{ip_forward} and \ref{ip_forward_agg} hold true
and assume that $G$ satisfies \ref{ip_G_invertible}.
Let $\phi$ and $\psi$
satisfy hypotheses \ref{ip-psiphi} Then the following estimate holds true:
\begin{equation}\label{stimabismut}
\vert Z^{t,x}_t\vert\leq C (T-t)^{-1/2},
\end{equation}
where $C$ depends on $t,\;T,\;A,\;F,\;\Vert\phi\Vert_\infty$.
\end{proposition}
\dim
We start by taking $\phi$ and $\psi$ differentiable with respect to their arguments.
Notice that the backward equation is differentiable e.g. following the results in \cite{BriFu}.

\noindent The proof follows in part the proof of Theorem 3.3 in \cite{Ri}, and we give it adequated to our setting.
By differentiating the forward-backward system in (\ref{fbsde}) with respect to the initial condition
$x$, we get
\begin{equation}\label{fbsde-differenziato}
    \left\{\begin{array}{l}\dis d\nabla_xX_\tau ^{t,x}=
A\nabla_xX_\tau ^{t,x}d\tau+\nabla_x F(\tau, X_\tau^{t,x})\nabla_xX_\tau^{t,x},\quad \tau\in
[t,T]\subset [0,T],
\\\dis
\nabla_xX_t^{t,x}=I,
\\\dis
 d\nabla_x Y_\tau^{t,x}=-\nabla_x\psi(\tau,X_\tau^{t,x},Y_\tau^{t,x},Z_\tau^{t,x})\nabla_xX_t^{t,x}\;d\tau
-\nabla_y\psi(\tau,X_\tau^{t,x},Y_\tau^{t,x},Z_\tau^{t,x})\nabla^xY_t^{t,x}\;d\tau\\
\qquad\qquad-\nabla_z\psi(\tau,X_\tau^{t,x},Y_\tau^{t,x},Z_\tau^{t,x})\nabla_xZ_t^{t,x}\;d\tau+
\nabla_xZ_\tau^{t,x}\;dW_\tau,
  \\\dis
 \nabla_x Y_T^{t,x}=\nabla_x\phi(X_T^{t,x})\nabla_xX_T^{t,x}.
\end{array}\right.
\end{equation}
By the Girsanov Theorem, there exists a probability measure $\Q$ such that the process $W^\Q$
given by
$$
W^\Q_\tau:=W_\tau-\int_t^\tau\nabla_z\psi(s,X_s^{t,x},Y_s^{t,x},Z_s^{t,x})\,ds
$$
is a cylindrical Wiener process. 
With respect to $\Q$ we can write the solution of the BSDE in (\ref{fbsde-differenziato})
as
\begin{align}\label{fbsde-differenziato-sol}
 \nabla_x Y^{t,x}_\tau&=\E_\tau^{\Q}\left[ e^{\int_\tau^T
\nabla_y \psi(s,X_s^{t,x},Y^{t,x}_s,Z^{t,x}_s)ds}\nabla\phi(X_T^{t,x})\nabla X_T^{t,x}\right.\\ \nonumber
&\left.+\int_\tau^Te^{\int_\tau^s
\nabla_y\psi(u,X_u^{t,x},Y^{t,x}_u,Z^{t,x}_u)du}\
\nabla_x\psi(s,X_s^{t,x},Y^{t,x}_s,Z^{t,x}_s)\nabla X_s^{t,x}ds\right].
\end{align}
Notice that with this change of measure, the equation satisfied by
the process $\nabla_xX^{t,x}$ does not change because the equation satisfied by $\nabla_xX^{t,x}$ is deterministic, the noise does not enter in this equation, and by proposition \ref{prop-ex-forward}
we get that $\nabla_x X^{t,x}$ is a bounded process. Moreover also notice that by our assumptions on
\ref{ip-psiphi}, $\psi$ is lipschitz continuous with respect to $x$ and $y$, and so  $\nabla_x\psi$ and $\nabla_y\psi$ are
bounded by $L_\psi$.
By setting
\begin{align}\label{def-F}
F^{t,x}_\tau=&e^{\int_t^\tau \nabla_y \psi(s,X_s^{t,x},Y^{t,x}_s,Z^{t,x}_s)ds}\nabla_x
Y_\tau^{t,x}\\
&+\int_t^\tau\nabla_x\psi(s,X_s^{t,x},Y^{t,x}_s,Z^{t,x}_s)\nabla X_s^{t,x}
e^{\int_t^s \nabla_y \psi(u,X_u^{t,x},Y^{t,x}_u,Z^{t,x}_u)du}\,ds,\nonumber
\end{align}
we get, by (\ref{fbsde-differenziato-sol}),
\[
 F^{t,x}_\tau=F_T^{t,x}-\int_\tau^Te^{\int_t^s
\nabla_y\psi(u,X_u^{t,x},Y^{t,x}_u,Z^{t,x}_u)du}\nabla Z_s^{t,x}\, dW^\Q_s,
\]
so $F^{t,x}$ is a $Q$-martingale and consequently $\left(F^{t,x}\right)^2$ is a $Q$-submartingale, so the following inequalities holds true
\begin{equation}\label{stima1}
 \E^\Q \int_\tau^T\vert F_s^{t,x}\vert^2\,ds\geq (T-t) \vert F^{t,x}_t\vert^2=
(T-t)\vert \nabla_xY^{t,x}_t\vert^2=(T-t) \vert Z^{t,x}_tG^{-1}(t)\vert^2;
\end{equation}
Moreover by (\ref{def-F}) we can deduce an expression for $\E^\Q \int_t^T\vert F_s^{t,x}\vert^2\,ds$ by noticing that
\begin{align*}
 &\E^\Q \int_t^T\vert F_s^{t,x}\vert^2\,ds=  \E^\Q \int_\tau^T\E^\Q_s\vert F_s^{t,x}\vert^2\,ds
\end{align*}
Taking into account that $Z^{t,x}$ is $\Q$-square integrable, that $\vert \nabla_x \psi\vert\leq L_\psi$ and
$\vert \nabla_y\psi \vert\leq L_\psi$, and that $\nabla_x X^{t,x}$ is a bounded process,
 we get that
\[
 \E^\Q \int_t^T\vert F_s^{t,x}\vert^2\,ds\leq C
\]
where $C$ depends on the coefficient of the forward backward system but does not depend on $\nabla\phi$.
Putting together this estimate with the estimate in (\ref{stima1}), and taking into account that $G$ has a bounded inverse,
we finally get
\begin{equation}\label{stima2}
\vert Z^{t,x}_t\vert\leq C (T-t)^{-1/2},
\end{equation}
where $C$ depends on the coefficient of the forward backward system but does not depend on $\nabla\phi$.

\noindent So estimate (\ref{stimabismut}) has been proved in the case of generator $\psi$
and final datum $\phi$ differentiable.
If they are not differentiable, we can approximate $\phi $ and $\psi$ with their inf-sup
convolutions $\phi_n$ and $\psi_n$ respectively, where $\phi_n$ is given by
\begin{equation}
\phi_{n}\left(  x\right)   =\sup_{x_2\in
H}\left\{  \inf_{x_1\in H}\left[  \phi\left(  x_1\right)  +\frac{n\left|
x_2-x_1\right|  _{H}^{2}}{2}\right]  -n\left|  x-x_2\right|  _{H}^{2}\right\}
  . \label{infsupconvphi}%
\end{equation}
For what concerns $\psi_n$, following e.g. the appendix in \cite{Mas-infor}, we define
\[
\bar\psi\left( t, x,y,z\right)=\dfrac{\psi\left( t, x,y,z\right)}{1+\vert z\vert^2}
\]
and 
\begin{align*}
\bar\psi_n\left( t, x,y,z\right)
 &=\sup_{x_1\in H,\,y_1\in\R,\,z_1 \in H}\Big\lbrace \inf_{x_2\in H,\,y_2\in\R,\,z_2 \in H}\Big[ 
   \bar\psi\left(t,  x_2,y_2,z_2\right)
\\ \nonumber 
 +&\dfrac{n\left(|
   x_1-x_2| _H^2+|
   y_1-y_2|+|
   z_1-z_2| _H^2\right)}{2}\Big] -\left(|  x-x_1|  _H^2
 -|  y-y_1|^2-|  z-z_1|  _H^2\right)\Big\rbrace.\\ \nonumber
\end{align*}
Finally we set
\begin{align}\label{infsupconvpsi}
\psi_n\left( t, x,y,z\right)=\bar\psi_n\left( t, x,y,z\right)\left(1+\vert z\vert^2\right)
\end{align}
By properties of inf-sup convolutions, see e.g. \cite{DP3}, we know that 
$\phi_n$ and $\psi_n$ are differentiable, and that they preserve the lipschitz constant, so that
\[
 \vert \nabla_x \psi_n \vert \leq L_\psi,\quad\vert \nabla_y \psi_n \vert \leq L_\psi.
\]
Moreover, also $\psi_n$ is differentiable, and by standard properties of the inf-sup convolutions of bounded functions, as well as properties of the inf-sup convolutions of functions with polynomial growth, see again
\cite{Mas-infor}, we have that
\[
 \sup_{x\in H,\,y\in\R,\,z \in H} \vert \dfrac{\psi_n(t,x,z)}{1+\vert z\vert^2}-
\dfrac{\psi(t,x,z)}{1+\vert z\vert^2}\vert \rightarrow 0 \quad \text{as }n\rightarrow\infty.
\]
We define $(Y^{n,t,x}, Z^{n,t,x})$ solutions of a BSDE like the one in the forward-backward system
(\ref{fbsde}) with generator $\psi_n$ and final datum $\phi_n$, namely
\begin{equation*}
\left\lbrace\begin{array}{l}
 -dY^{n,t,x}_\tau=\psi_n(\tau,X_\tau^{t,x},Y^{n,t,x}_\tau,Z^{n,t,x}_\tau)\,d\tau-Z^{n,t,x}_\tau dW_\tau\\
Y^{n,t,x}_T=\phi_n(X_T^{t,x}).
\end{array}\right.
\end{equation*}
It is a known result that 
\[
 \Vert \left( Y^{n,t,x}- Y^{t,x},Z^{n,t,x}- Z^{t,x} \right)\Vert_{\cals^2\times\calm^2}\rightarrow 0 \qquad \text{as }n\rightarrow\infty.
\]
Moreover for $Z^{n,t,x}$ we can prove an estimate like (\ref{stima2}), where $C$ does not depend on $n$.
Since, at least by taking a subsequence, $Z_{n_{k}}\rightarrow Z $ $\P$-almost surely, we finally conclude that (\ref{stimabismut})
holds true for $Z$, with $\phi$ and $\psi$ satisfying hypothesis \ref{ip-psiphi}, and the proposition is proved.
\qed

We conclude this section by proving some further integrability properties of the process $Z$, namely we prove that $Z\in\calm^p$, for any $p\geq 1$. The proof of this result is quickly given e.g. in \cite{BriFu}, section 4, inequality $(13)$, and the proof is based on Kobylanski transform introduced in \cite{Kob}.
We give here an alternative proof, and in the next proposition
\ref{prop-convp-bsdequadr} we also prove the stability of the solution
with respect to the final datum in $\cals^\infty\times \calm^p$-norm.
\begin{proposition}\label{teo-stimap-bsdequadr}Let $(X,Y,Z)$ be solution of the forward-backward system (\ref{fbsde}),
and assume that hypotheses \ref{ip_forward}, \ref{ip_forward_agg}, \ref{ip-psiphi}
hold true and assume that $G$ satisfies \ref{ip_G_invertible}.
Then, for all $p\geq 1$, the unique solution of the markovian BSDE in (\ref{fbsde}) is such that
\[
 \Vert Y\Vert _{\cals^\infty}+\Vert Z\Vert_{\calm^p}\leq C,
\]
where $C$ is a constant that may depend on $A,\,F,\,G,\, K_\psi,\,L_\psi,\,K_\phi, t,\,T$.
\end{proposition}
\dim  We start by proving that $Y$ is a bounded process. We notice that $\forall t\leq \tau\leq T$,
by the Markov property we have $X_\tau^{t,x}=X_\tau^{\tau,y}\arrowvert_{y=X_\tau^{t,x}}$ and
\[
 Y_\tau^{t,x}=Y_\tau^{t,X_\tau^{t,x}}=Y_\tau^{\tau,y}\arrowvert_{y=X_\tau^{t,x}}, \qquad Z_\tau^{t,x}=Z_\tau^{t,X_\tau^{t,x}}=Z_\tau^{\tau,y}\arrowvert_{y=X_\tau^{t,x}}.
\]
Moreover, by writing the BSDE in the forward-backward system (\ref{fbsde}) in integral form and with initial conditions given by $\tau$ and $y$ we get
\[
Y_\tau^{\tau,y} =\phi(X_T^{\tau,y})+\int_\tau^T\psi\left(r,X_r^{\tau,y},Y_r^{\tau,y},Z_r^{\tau,y}\right)
\,dr-\int_\tau^T  Z_r^{\tau, y}\,dW_r.
\]
By taking expectation, and by taking into account that by \ref{teo-ex-bsdequadr}, $(Y,Z)\in \cals^2\times
\calm^2$,  and noting that $Y_\tau^{\tau,y}$ is deterministic,  we get
\begin{equation*}
\vert  Y_\tau^{\tau,y}\vert \leq\E\vert\phi(X_T^{\tau,y})\vert
+\E\int_\tau^T\vert\psi\left(r,X_r^{\tau,y},Y_r^{\tau,y},Z_r^{\tau,y}\right)\vert
\,dr+\E\left(\int_\tau^T\vert Z^{t,x}_r\vert^2\,dr\right)^{1/2}
\end{equation*}
where $C$ is a constant that may depend on $A,\,F,\,G,\, K_\psi,\,L_\psi,\,K_\phi, t,\,T$ but not on $x$.
We have also
\[
 \vert Y_\tau^{t,x}\vert=\vert Y_\tau^{\tau,X_\tau^{t,x}} \vert\leq C,
\]
so that we have proved that $Y^{t,x}\in\cals^\infty$.
Now we have to prove that $Z^{t,x}\in\calm^p$, for every $p\geq 1$. We already know that
$Z^{t,x}\in\calm^2$, so it remains to prove that $Z^{t,x}\in\calm^p$, for $p> 2$. To this aim,
we notice that by applying It\^o's formula to $\vert Y^{t,x}\vert^2$ and by integrating over $[t,T-\delta]$,
for $\delta>0$ arbitrarly chosen,
we obtain
\begin{align*}
\vert &Y_t^{t,x}\vert^2 +\int_t^{T-\delta}\vert Z^{t,x}_s \vert^2 \,ds\\ \nonumber
&=\vert Y_{T-\delta}^{t,x}\vert^2+2\int_t^{T-\delta}Y_s^{t,x}\psi\left(s,X_s^{t,x},Y_s^{t,x},Z_s^{t,x}\right)
\,ds- \int_t^{T-\delta}  2Y_s^{t,x}Z_s^{t,x}\,dWs.
\end{align*}
Raising to the power $p/2$ and taking expectation we get
\begin{align}\label{stima1p}
\vert& Y_t^{t,x}\vert^p +\E\left(\int_t^{T-\delta}\vert Z^{t,x}_s \vert^2 \,ds\right)^{p/2}\\ \nonumber
&\leq C_p\left[\E\vert Y_{T-\delta}^{t,x}\vert^p+\E\vert\int_t^{T-\delta}Y_s^{t,x}\psi\left(s,X_s^{t,x},Y_s^{t,x},Z_s^{t,x}\right)
\,ds\vert^{p/2}+ \E\vert\int_t^{T-\delta}  Y_s^{t,x}Z_s^{t,x}\,dW_s\vert^{p/2}\right].
\end{align}
where $C$ is a constant that may depend on $p$.
We start by estimating $\E\vert\int_t^{T-\delta}  Y_s^{t,x}Z_s^{t,x}\,dW_s\vert^{p/2}$: by the Burkholder-Davies-Gundy inequality
and since $Y^{t,x}$ is a uniformly bounded process, we get
\begin{align}\label{stimaYZp}
 \E&\vert\int_t^{T-\delta}  Y_s^{t,x}Z_s^{t,x}\,dW_s\vert^{p/2}\leq
\E\left(\int_t^{T-\delta} \vert Y_s^{t,x}Z_s^{t,x}\vert^2\,ds\right)^{p/4}\\ \nonumber
&\leq \E C^{p/2}\left(\int_t^{T-\delta} \vert Z_s^{t,x}\vert^2\,ds\right)^{p/4}\leq a_1 C^{p}+a_2
\E\left(\int_t^{T-\delta} \vert Z_s^{t,x}\vert^2\,ds\right)^{p/2}
\end{align}
where in the last passage we have applied Young inequality and $a_1*a_2=1/4$, and we will choose later $a_1$
and $a_2$ such that $a_2$ is sufficiently small.

\noindent Next we estimate $\E\vert\int_t^{T-\delta}Y_s^{t,x}\psi\left(s,X_s^{t,x},Y_s^{t,x},Z_s^{t,x}\right)
\,ds\vert^{p/2}$: by hypothesis \ref{ip-psiphi} it follows
\[
 \vert\psi(t,x,y,z)\vert\leq C\left(1+\vert y\vert+\vert z\vert^2  \right),
\]
so we get
\begin{align}\label{stimaYZp1}
 \E&\vert\int_t^{T-\delta}Y_s^{t,x}\psi\left(s,X_s^{t,x},Y_s^{t,x},Z_s^{t,x}\right)
\,ds\vert^{p/2}\\ \nonumber
&\leq C\E\left(\int_t^{T-\delta}\left( 1+ \vert Y_s^{t,x}\vert + \vert Y_s^{t,x}\vert^2\right)
\,ds\right)^{p/2}+ C\E\left(\int_t^{T-\delta}\vert Y_s^{t,x}\vert \vert Z_s^{t,x}\vert^2
\,ds\right)^{p/2}\\ \nonumber
&\leq C+C\E\left(\int_t^{T-\delta}\vert Y_s^{t,x}\vert \vert Z_s^{t,x}\vert^2
\,ds\right)^{p/2},
\end{align}
where the last passage follows since $Y^{t,x}$ is a bounded process. We have to estimate the last integral
in (\ref{stimaYZp1}): 
by estimate (\ref{stimabismut}), we get that on the interval $[t, T-\delta]$
$Z^{t,x}$ satisfies
\[
 \vert Z_t^{\tau,x}\vert\leq C(T-\tau)^{-1/2}\leq C\delta^{-1/2}\qquad\text{for }\tau\in[t,T-\delta].
\]
So
\begin{align}\label{stimaYZp1bis}
 \E&\left(\int_t^{T-\delta}\vert Y_s^{t,x}\vert \vert Z_s^{t,x}\vert^2
\,ds\right)^{p/2}\\ \nonumber
&\leq \delta^{-(\delta/2)( p/2)}\E\left(\int_t^{T-\delta}\vert Y_s^{t,x}\vert \vert Z_s^{t,x}\vert^{2-\delta}
\,ds\right)^{p/2}\\ \nonumber
& \leq C \delta^{-(\delta p)/4}\E\left(\int_t^{T-\delta} \vert Z_s^{t,x}\vert^{2-\delta}
\,ds\right)^{p/2}.
\end{align}
By Young inequality $\vert Z_s^{t,x}\vert^{2-\delta}\leq \bar a _1 (\delta/2) 1^{2/\delta}+\bar a _2 (2-\delta)/2 \vert Z_s^{t,x}\vert^{2} $ with  $\bar a_1*\bar a_2=1$, and we will choose later $\bar a_1$
and $\bar a_2$ such that $\bar a_2$ is sufficiently small.
So, going on with estimate (\ref{stimaYZp1}), also noting that $(2-\delta)/2<1$ we get
\begin{align}\label{stimaYZp2}
 \E&\left(\int_t^{T-\delta}\vert Y_s^{t,x}\vert \vert Z_s^{t,x}\vert^2
\,ds\right)^{p/2}\\ \nonumber
&\leq C_p\delta^{-(\delta p)/4}\bar a_1^{p/2} \left(\delta/2\right)^{p/2}
+C_p\bar a_2^{p/2} \delta^{-(\delta p)/4}\E\left(\int_t^{T-\delta} \vert Z_s^{t,x}\vert^{2}
\,ds\right)^{p/2}.
\end{align}
Coming back to (\ref{stima1p}), by estimate (\ref{stimaYZp}) and (\ref{stimaYZp2}), and since for $\delta $ small $\delta^{-(\delta p)/4}$ is uniformly bounded by $1$, we get
\begin{align*}
\vert& Y_t^{t,x}\vert^p +\E\left(\int_t^{T-\delta}\vert Z^{t,x}_s \vert^2 \,ds\right)^{p/2}\\ \nonumber
&\leq C\left(1+a_1+ \delta^{-(\delta p)/4}\bar a_1^{p/2} \left(\delta/2\right)^{p/2}\right)
+C\left(a_2+ \bar a_2^{p/2} \right)\E\left(\int_t^{T-\delta} \vert Z_s^{t,x}\vert^{2}
\,ds\right)^{p/2}.
\end{align*}
By choosing $a_1,\,a_2,\,\bar a_1,\,\bar a_2$ such that $C\left(a_2+ \bar a_2^{p/2} \right)\leq1/2$
and since $\delta^{-(\delta p)/4}\left(\delta/2\right)^{p/2}$ is uniformly bounded in $\delta$,
we finally get that
\[
 \E\left(\int_t^{T-\delta}\vert Z^{t,x}_s \vert^2 \,ds\right)^{p/2}\leq C,
\]
where $C$ does not depend on $\delta$, so letting $\delta$ go to $0$, by the monotone convergence theorem we finally get
\[
 \E\left(\int_t^{T}\vert Z^{t,x}_s \vert^2 \,ds\right)^{p/2}\leq C
\]
and the proof is concluded.
\qed

Finally we prove a stability result with respect to approximation of the final datum
$\phi$ for the $\calm^p$-norm of $Z$ in the next propositon, and in the succeeding one
we will prove a stability result with respect to approximation of the generator
$\psi$ for the $\calm^p$-norm of $Z$.
\begin{proposition}\label{prop-convp-bsdequadr}Let $(X,Y,Z)$ be solution of the forward-backward system (\ref{fbsde}),
and assume that hypotheses \ref{ip_forward} and \ref{ip-psiphi} hold
true and let $(Y^n,Z^n)$ be solution of the BSDE in the forward-backward system (\ref{fbsde})
with final datum equal to $\phi_n$ in the place of $\phi$, and such that $\forall n\geq 1$
$\phi_n$ satisfies \ref{ip-psiphi} and $\Vert \phi_n\Vert_\infty\leq\Vert\phi\Vert_\infty$; and assume
that $\Vert \phi_n-\phi\Vert_\infty\rightarrow 0$
as $n\rightarrow \infty$.
Then, for all $p\geq 1$, the unique solution of the markovian BSDE in (\ref{fbsde}) is such that
\[
 \Vert Y-Y^n\Vert _{\cals^p}+\Vert Z-Z^n\Vert_{\calm^p}\rightarrow 0 \qquad \text{as }n\rightarrow \infty.
\]
\end{proposition}
\dim It is well known that under these assumptions 
\[
 \Vert Y-Y^n\Vert _{\cals^2}+\Vert Z-Z^n\Vert_{\calm^2}\rightarrow 0 \qquad \text{as }n\rightarrow \infty.
\]
and since by theorem \ref{teo-stimap-bsdequadr} $Y$ is a bounded process and the sequence $(Y^n)_n$
is a sequence uniformly bounded with respect to $n$,
it immediately follows that $\Vert Y-Y^n\Vert _{\cals^p}\rightarrow 0$ as $n\rightarrow \infty$, for any $p\geq 1$. We also notice that the assumptions on $\phi_n$ in the proposition are satisfied by the
inf-sup convolutions of $\phi$.

Now we have to prove the convergence of  $Z^{n,t,x}$ to $Z^{t,x}$ in $\calm^p$. Similarly
to the proof of theorem \ref{teo-stimap-bsdequadr}
we apply It\^o's formula to $\vert Y^{n,t,x}-Y^{t,x}\vert^2$, we integrate over $[t,T-\delta]$,
for $\delta>0$ arbitrarly chosen, we raise to the power $p/2$ and taking expectation we get
\begin{align}\label{stima1pn}
\vert& Y_t^{n,t,x}- Y_t^{t,x}\vert^p +\E\left(\int_t^{T-\delta}\vert  Z_t^{n,t,x}- Z^{t,x}_s \vert^2 \,ds\right)^{p/2}\\ \nonumber
&\leq C_p\left[\E\vert  Y_{T-\delta}^{n,t,x} -Y_{T-\delta}^{t,x}\vert^p+ \E\vert\int_t^{T-\delta} \left(Y_s^{n,t,x} -Y_s^{t,x}\right)
\left(Z_s^{n,t,x}-Z_s^{t,x}\right)\,dW_s\vert^{p/2}\right.\\ \nonumber
&\left.+\E\vert\int_t^{T-\delta}\left( Y_s^{n,t,x}-Y_s^{t,x}\right)\left(\psi\left(s,X_s^{t,x},Y_s^{n,t,x},Z_s^{n,t,x}\right)
-\psi\left(s,X_s^{t,x},Y_s^{t,x},Z_s^{t,x}\right)\right)
\,ds\vert^{p/2}\right].
\end{align}
where $C_p$ is a constant that may depend on $p$.
The estimate of
$$\E\vert\int_t^{T-\delta} \left(Y_s^{n,t,x} Y_s^{t,x}\right)
\left(Z_s^{n,t,x}-Z_s^{t,x}\right)\,dW_s\vert^{p/2}
$$ can be performed exactly as the estimate of $\E\vert\int_t^{T-\delta}  Y_s^{t,x}Z_s^{t,x}\,dW_s\vert^{p/2}$ in the proof of the previous theorem \ref{teo-stimap-bsdequadr}, see (\ref{stimaYZp}), arriving at
\begin{align*}
 \E&\vert\int_t^{T-\delta}  \left(Y_s^{n,t,x}-Y_s^{t,x}\right)\left(Z_s^{n,t,x}-Z_s^{t,x}\right)\,dW_s\vert^{p/2}
\leq\E\left(\int_t^{T-\delta} \vert Y_s^{t,x}Z_s^{t,x}\vert^2\,ds\right)^{p/4}\\ \nonumber
&\leq a_2\E\left(\int_t^{T-\delta} \vert Z_s^{t,x}-Z_s^{n,t,x}\vert^2\,ds\right)^{p/2}+a_1 
\E\sup_{s\in[t,T]}\vert Y_s^{t,x}-Y_s^{n,t,x}\vert^p
\E\left(\int_t^{T-\delta} \vert Z_s^{t,x}\vert^2\,ds\right)^{p/2}
\end{align*}

\noindent Now we estimate, using also hypothesis \ref{ip-psiphi},
\begin{align}\label{stimaYZp1conv}
 \E&\vert\int_t^{T-\delta}\left(Y_s^{n,t,x}-Y_s^{t,x}\right)
\left(\psi\left(s,X_s^{t,x},Y_s^{n,t,x},Z_s^{n,t,x}\right)-\psi\left(s,X_s^{t,x},Y_s^{t,x},Z_s^{t,x}\right)\right)
\,ds\vert^{p/2}\\ \nonumber
&\leq C\E\left(\int_t^{T-\delta} \vert Y_s^{n,t,x}- Y_s^{t,x}\vert^2
\,ds\right)^{p/2}\\ \nonumber
&+ C\E\left(\int_t^{T-\delta}\vert Y_s^{n,t,x}-Y_s^{t,x}\vert \vert Z_s^{t,x}-Z_s^{n,t,x}\vert
\left(1+\vert Z_s^{t,x}\vert+\vert Z_s^{n,t,x}\vert  \right)
\,ds\right)^{p/2}\\ \nonumber
&\leq C+Cc_p\E\left(a_1\int_t^{T-\delta}\vert Y_s^{n,t,x}-Y_s^{t,x}\vert^2 
\left(1+\vert Z_s^{t,x}\vert+\vert Z_s^{n,t,x}\vert  \right)^2
\,ds\right)^{p/2}\\ \nonumber
&+Cc_p\E\left(a_2\int_t^{T-\delta}\vert Z_s^{t,x}-Z_s^{n,t,x}\vert^2
\,ds\right)^{p/2},
\end{align}
where we have applied Young inequality, and $a_1$ and $a_2$ have been chosen
such that  $a_1*a_2=1/4$, and $C c_p a_2^{p/2}<1/2$. We have to estimate the first integral
in the last passage of (\ref{stimaYZp1conv}): 
by (\ref{stimabismut}), we get that on the interval $[t, T-\delta]$,
$Z^{n,t,x}$ and $Z^{t,x}$ satisfy
\[
 \vert Z^{n,t,x}_\tau\vert+\vert Z_\tau^{t,x}\vert\leq C(T-\tau)^{-1/2}\leq C\delta^{-1/2}\qquad\text{for }\tau\in[t,T-\delta],
\]
with $C$ that does not depend on $n$.
So
\begin{align}\label{stimaYZp1convbis}
 \E&\left(a_1\int_t^{T-\delta}\vert Y_s^{n,t,x}-Y_s^{t,x}\vert^2 
\left(\vert Z_s^{t,x}\vert+\vert Z_s^{n,t,x}\vert^2  \right)
\,ds\right)^{p/2}\\ \nonumber
 &\leq \delta^{-(\delta/2)( p/2)}\E\left(\int_t^{T-\delta}\vert Y_s^{n,t,x}-Y_s^{t,x}\vert^2 \left(
 \vert Z_s^{n,t,x}\vert^{2-\delta}+\vert Z_s^{t,x}\vert^{2-\delta}\right)
 \,ds\right)^{p/2}\\ \nonumber
 & \leq  \delta^{-(\delta p)/4}\E\left[\left(\int_t^{T-\delta}\vert Y_s^{n,t,x}-Y_s^{t,x}\vert^{4/\delta}\,ds\right)^{(p\delta)/4 }\right.\\
&\left.\left( \int_t^{T-\delta}\left(\vert Z_s^{n,t,x}\vert^{2-\delta}+\vert Z_s^{t,x}\vert^{2-\delta}\right)^{2/(2-\delta)}
 \,ds\right)^{p(2-\delta)/4}\right]\\ \nonumber
 & \leq  \delta^{-(\delta p)/4}\left(\E\left(\int_t^{T-\delta}\vert Y_s^{n,t,x}-Y_s^{t,x}\vert^{4/\delta}\,ds\right)^{(p\delta)/2}\right)^{1/2} \\ \nonumber
&\qquad\left(\E\left(
 \int_t^{T-\delta}\left(\vert Z_s^{n,t,x}\vert^{2}+\vert Z_s^{t,x}\vert^{2}\right)
 \,ds\right)^{p(2-\delta)/2}\right)^{1/2}\\ \nonumber
 &\leq C\E\sup_{s\in[t,T]}\vert Y_s^{n,t,x}-Y_s^{t,x}\vert^p.
\end{align}
We finally get that
\[
 \E\left(\int_t^{T-\delta}\vert Z^{n,t,x}_s-Z^{t,x}_s \vert^2 \,ds\right)^{p/2}\leq C\E\sup_{s\in[t,T]}\vert Y_s^{n,t,x}-Y_s^{t,x}\vert^p,
\]
where $C$ does not depend on $\delta$, so letting $\delta$ go to $0$, by the monotone convergence theorem we finally get
\[
 \E\left(\int_t^{T}\vert Z^{t,x}_s -Z^{n,t,x}_s\vert^2 \,ds\right)^{p/2}\leq C\E\sup_{s\in[t,T]}\vert Y_s^{n,t,x}-Y_s^{t,x}\vert^p,
\]
and the proof is concluded.
\qed
\begin{proposition}\label{prop-psiconvp-bsdequadr}Let $(X,Y,Z)$ be solution of the forward-backward system (\ref{fbsde}),
and assume that hypotheses \ref{ip_forward} and \ref{ip-psiphi} hold
true and let $(Y^n,Z^n)$ be solution of the BSDE in the forward-backward system (\ref{fbsde})
with generator  $\psi_n$ in the place of $\psi$, where $\psi_n$ is defined in (\ref{infsupconvpsi}). Namely,
$(Y^n,Z^n)$ solve the following BSDE:
\begin{equation}\label{bsde-n-psi}
    \left\{\begin{array}{l}\dis
 dY^{n,t,x}_\tau=-\psi_n(\tau,X^{t,x}_\tau,Y^{n,t,x}_\tau,Z^{n,t,x}_\tau)\;d\tau+Z^{n,t,x}_\tau\;dW_\tau,
  \\\dis
  Y^{n,t,x}_T=\phi(X_T),
\end{array}\right.
\end{equation}
Then, for all $p\geq 1$, the unique solution of the markovian BSDE in (\ref{fbsde}) is such that
\[
 \Vert Y-Y^n\Vert _{\cals^p}+\Vert Z-Z^n\Vert_{\calm^p}\rightarrow 0 \qquad \text{as }n\rightarrow \infty.
\]
\end{proposition}
\dim As in the proof of the previous proposition, it is well known that under these assumptions 
\[
 \Vert Y-Y^n\Vert _{\cals^2}+\Vert Z-Z^n\Vert_{\calm^2}\rightarrow 0 \qquad \text{as }n\rightarrow \infty,
\]
and consequently $\Vert Y-Y^n\Vert _{\cals^p}\rightarrow 0$ as $n\rightarrow \infty$, for any $p\geq 1$.
Now we have to prove the convergence of  $Z^{n,t,x}$ to $Z^{t,x}$ in $\calm^p$. Similarly
to the proof of theorem \ref{teo-stimap-bsdequadr}
we apply It\^o's formula to $\vert Y^{n,t,x}-Y^{t,x}\vert^2$, we integrate over $[t,T-\delta]$,
for $\delta>0$ arbitrarly chosen, we raise to the power $p/2$ and taking expectation we get
\begin{align*}
&Y_t^{n,t,x}- Y_t^{t,x}\vert^p +\E\left(\int_t^{T-\delta}\vert  Z_t^{n,t,x}- Z^{t,x}_s \vert^2 \,ds\right)^{p/2}\\ \nonumber
&\leq C_p\E\vert  Y_{T-\delta}^{n,t,x} -Y_{T-\delta}^{t,x}\vert^p+ \E\vert\int_t^{T-\delta} \left(Y_s^{n,t,x} -Y_s^{t,x}\right)
\left(Z_s^{n,t,x}-Z_s^{t,x}\right)\,dW_s\vert^{p/2}\\ \nonumber
&+\E\vert\int_t^{T-\delta}\left( Y_s^{n,t,x}-Y_s^{t,x}\right)\left(\psi_n\left(s,X_s^{t,x},Y_s^{n,t,x},Z_s^{n,t,x}\right)
-\psi\left(s,X_s^{t,x},Y_s^{t,x},Z_s^{t,x}\right)\right)
\,ds\vert^{p/2}.
\end{align*}
where $C$ is a constant that may depend on $p$.
The estimate of 
$$\E\vert\int_t^{T-\delta} \left(Y_s^{n,t,x} Y_s^{t,x}\right)
\left(Z_s^{n,t,x}-Z_s^{t,x}\right)\,dW_s\vert^{p/2}
$$ can be performed exactly as the estimate of $\E\vert\int_t^{T-\delta}  Y_s^{t,x}Z_s^{t,x}\,dW_s\vert^{p/2}$ in the proof of the previous theorem \ref{teo-stimap-bsdequadr}.

\noindent Now we estimate, using also hypothesis \ref{ip-psiphi} and property of the inf-sup convolutions,,
\begin{align}\label{stimaYZp1convpsi}
 \E&\vert\int_t^{T-\delta}\left(Y_s^{n,t,x}-Y_s^{t,x}\right)
\left(\psi_n\left(s,X_s^{t,x},Y_s^{n,t,x},Z_s^{n,t,x}\right)-\psi\left(s,X_s^{t,x},Y_s^{t,x},Z_s^{t,x}\right)\right)
\,ds\vert^{p/2}\\ \nonumber
&\leq \E\vert\int_t^{T-\delta}\left(Y_s^{n,t,x}-Y_s^{t,x}\right)
\left(\psi_n\left(s,X_s^{t,x},Y_s^{n,t,x},Z_s^{n,t,x}\right)-\psi\left(s,X_s^{t,x},Y_s^{n,t,x},Z_s^{n,t,x}\right)\right)
\,ds\vert^{p/2}\\ \nonumber
&+\E\vert\int_t^{T-\delta}\left(Y_s^{n,t,x}-Y_s^{t,x}\right)
\left(\psi\left(s,X_s^{t,x},Y_s^{n,t,x},Z_s^{n,t,x}\right)-\psi\left(s,X_s^{t,x},Y_s^{t,x},Z_s^{t,x}\right)\right)
\,ds\vert^{p/2}\\ \nonumber
&\leq \E\vert\int_t^{T-\delta}\left(Y_s^{n,t,x}-Y_s^{t,x}\right)
\left(\psi_n\left(s,X_s^{t,x},Y_s^{n,t,x},Z_s^{n,t,x}\right)-\psi\left(s,X_s^{t,x},Y_s^{n,t,x},Z_s^{n,t,x}\right)\right)
\,ds\vert^{p/2}\\ \nonumber
&+C\E\left(\int_t^{T-\delta} \vert Y_s^{n,t,x}- Y_s^{t,x}\vert^2
\,ds\right)^{p/2}\\ \nonumber
&+ C\E\left(\int_t^{T-\delta}\vert Y_s^{n,t,x}-Y_s^{t,x}\vert \vert Z_s^{t,x}-Z_s^{n,t,x}\vert
\left(1+\vert Z_s^{t,x}\vert+\vert Z_s^{n,t,x}\vert  \right)
\,ds\right)^{p/2}\\ \nonumber
&\leq\left(
 \sup_{x\in H,\,y\in\R,\,z \in H} \vert \dfrac{\psi_n(t,x,y,z)}{1+\vert z\vert^2}-
\dfrac{\psi(t,x,y,z)}{1+\vert z\vert^2}\vert 
\right)^{p/2}\E\left(\int_t^{T-\delta}\vert Y_s^{n,t,x}-Y_s^{t,x}\vert\left( 1+\vert Z^{n,t,x}_s\vert^2\right)
\,ds\right)^{p/2}\\ \nonumber
& +Cc_p\E\left(a_2\int_t^{T-\delta}\vert Z_s^{t,x}-Z_s^{n,t,x}\vert^2
\,ds\right)^{p/2}\\ \nonumber
&+Cc_p\E\left(a_1\int_t^{T-\delta}\vert Y_s^{n,t,x}-Y_s^{t,x}\vert^2 
\left(1+\vert Z_s^{t,x}\vert+\vert Z_s^{n,t,x}\vert  \right)^2
\,ds\right)^{p/2}\\ \nonumber
&\leq\left(
 \sup_{x\in H,\,y\in\R,\,z \in H} \vert \dfrac{\psi_n(t,x,y,z)}{1+\vert z\vert^2}-
\dfrac{\psi(t,x,y,z)}{1+\vert z\vert^2}\vert 
\right)^{p/2}C^{p/2}\E\left(\int_t^{T-\delta}\vert\left( 1+\vert Z^{n,t,x}_s\vert^2\right)
\,ds\right)^{p/2}\\ \nonumber
& +Cc_p\E\left(a_2\int_t^{T-\delta}\vert Z_s^{t,x}-Z_s^{n,t,x}\vert^2
\,ds\right)^{p/2}\\ \nonumber
&+Cc_p\E\left(a_1\int_t^{T-\delta}\vert Y_s^{n,t,x}-Y_s^{t,x}\vert^2 
\left(1+\vert Z_s^{t,x}\vert+\vert Z_s^{n,t,x}\vert \right)^2
\,ds\right)^{p/2},
\end{align}
where we have applied Young inequality, and $a_1$ and $a_2$ are choosen so that
 $a_1*a_2=1/2$, and $C c_p a_2^{p/2}<1/2$. We notice that
\[\E\vert\int_t^{T-\delta}\vert Y_s^{n,t,x}-Y_s^{t,x}\vert
\,ds\vert^{p/2}\rightarrow 0 \qquad \text{as }n\rightarrow \infty, \text{ and } \E\int_t^T\vert Z^{n,t,x}_s\vert^2\,ds\leq C,
\]
with $C$ independent on $n$,
and also that
\[
 \sup_{x\in H,\,y\in\R,\,z \in H} \vert \dfrac{\psi_n(t,x,z)}{1+\vert z\vert^2}-
\dfrac{\psi(t,x,z)}{1+\vert z\vert^2}\vert \rightarrow 0 \quad \text{as }n\rightarrow\infty.
\]
The other terms in the last passage of (\ref{stimaYZp1convpsi}) can be estimated as in (\ref{stimaYZp1conv}), finally arriving at
\[
 \E\left(\int_t^{T}\vert Z^{t,x}_s -Z^{n,t,x}_s\vert^2 \,ds\right)^{p/2}\leq C\E\sup_{s\in[t,T]}\vert Y_s^{n,t,x}-Y_s^{t,x}\vert^p,
\]
and the proof is concluded.
\qed

\subsection{The Bismut-Elworthy formula in the case of lipschitz generator}
\label{sez-Bismut-lip}
In this section we briefly recall the nonlinear version of the Bismut-Elworthy
formula proved in \cite{futeBismut}. To this aim, and also to prove the
Bismut-Elworthy formula in the quadratic case, which is the core of the paper,
we assume further that the operators $G(t)$
are boundedly invertible, as required in hypothesis \ref{ip_G_invertible}.

\noindent For $0\leq t\leq s\leq T$ and $h\in H$ we define the real valued random variables
\begin{equation}\label{def-U}
 U^{h,t,x}_s:=\dfrac{1}{s-t}\int_t^s \<G^{-1}(r,X^{t,x}_r)\nabla_xX_r^{t,x}h,dW_r \>
\end{equation}
We are ready to recall the Bismut-Elworthy formula proved in \cite{futeBismut}.
\begin{theorem}\label{teoBismutLip}
 Assume that hypotheses \ref{ip_forward}, \ref{ip_forward_agg}
and \ref{ip_G_invertible} hold true and let $\phi$ and $\psi$ in (\ref{fbsde}) be measurable, moreover for every fixed $t\in[0,T]$ the map $\psi(t,\cdot,\cdot,\cdot):H\times\R\times\Xi\rightarrow\R$
is continuous, and $\phi$ is continuous. Finally there exist nonegative constants
$L_\psi,\,K_\psi,\,K_\phi,\,\mu$
such that
\begin{align*}
&\vert \psi(t,x,y_1,z_1)- \psi(t,x,y_2,z_2)\vert
\leq L_\psi\left(\vert y_1-y_2\vert+\vert z_1-z_2\vert\right)
,\\
&\vert \psi(t,x,0,0)\vert\leq K_\psi\left(1+\vert x\vert^\mu\right), \qquad 
\vert \phi(x)\vert\leq K_\phi\left(1+\vert x\vert^\mu\right),
\end{align*}
for every $t\in[0,T]$, $x\in H$, $y_1,y_2\in\R$ and $z_1,z_2\in \Xi$
Then for $0\leq t\leq s\leq T$, $x,h\in H$
\begin{equation}\label{Bismut-lip}
 \E\left[ \nabla_x\,Y^{t,x}_sh \right]=
\E\int_s^T\psi\left(r,X_r^{t,x},Y_r^{t,x},Z_r^{t,x}\right)U^{h,t,x}_r\,dr
+\E\left[  \phi(X_T^{t,x}U^{h,x}_T\right]
\end{equation}
where 
\begin{equation}\label{bound-U}
 \left(\E\vert U^{h,t,x}_s \vert^2\right)^{1/2}\leq
C\left(s-t \right)^{-1/2}\vert h\vert.
\end{equation}
\end{theorem}
\dim
The proof is given in \cite{futeBismut}, lemma 3.8 and theorem 3.10.
\qed

In the following we will need also to generalize (\ref{bound-U}) from $q=2$ to any $q\geq 1$.
\begin{lemma}\label{lemma:boundU-q}
 Assume that hypotheses \ref{ip_forward}, \ref{ip_forward_agg}
and \ref{ip_G_invertible} hold true: for any $q\geq 1$,
\begin{equation}\label{bound-Uq}
 \left(\E\vert U^{h,t,x}_s \vert^q\right)^{1/q}\leq
C\left(s-t \right)^{-1/2}\vert h\vert, 
\end{equation}
and also
\begin{equation}
 \left(\E\sup_{s\in[\frac{t+T}{2},T]}\vert U^{h,t,x}_s \vert^q\right)^{1/q}\leq C\dfrac{1}{(T-t)^{1/2}}.\label{bound-Uq-sup}
\end{equation}
\end{lemma}
\dim We compute
\begin{align*}
 \E\vert U^{h,t,x}_s \vert^q =&
\E\vert \dfrac{1}{s-t}\int_t^s \<G^{-1}(r,X^{t,x}_r)\nabla_xX_r^{t,x}h,dW_r \>\vert^q\\ \nonumber
&\leq  \dfrac{1}{(s-t)^q}\E \left(\int_t^s \vert\<G^{-1}(r,X^{t,x}_r)\nabla_xX_r^{t,x}h\vert^2\,dr \>\right)^{q/2}\\ \nonumber
& \leq \dfrac{1}{(s-t)^q}\E \left(\int_t^s \vert\<G^{-1}(r,X^{t,x}_r)\nabla_xX_r^{t,x}h\vert^2\,dr \>\right)^{q/2}\\ \nonumber
& \leq \dfrac{1}{(s-t)^q}C (s-t)^{q/2}=C \dfrac{1}{(s-t)^{q/2}}.
\end{align*}
and also
\begin{align*}
 \E\sup_{s\in[\frac{t+T}{2},T]}\vert U^{h,t,x}_s \vert^q \leq&
\E\vert \dfrac{1}{(T-t)/2}\sup_{s\in[\frac{t+T}{2},T]}\int_t^s \<G^{-1}(r,X^{t,x}_r)\nabla_xX_r^{t,x}h,dW_r \>\vert^q\\ \nonumber
&\leq C  \dfrac{1}{(T-t)^q}\E \left(\int_t^T \vert\<G^{-1}(r,X^{t,x}_r)\nabla_xX_r^{t,x}h\vert^2\,dr \>\right)^{q/2}\\ \nonumber
& \leq C \dfrac{1}{(T-t)^{q/2}},
\end{align*}
which leads to 
\begin{equation*}
\left(\E\sup_{s\in[\frac{t+T}{2},T]}\vert U^{h,t,x}_s \vert^q\right)^{1/q}\leq C\dfrac{1}{(T-t)^{1/2}} 
\end{equation*}
\qed

\section{The Bismut-Elworthy formula in the quadratic case}
\label{sez-Bismut-quad}

This section is the core of the paper. We will work with a generator $\psi$ with quadratic growth with respect to $z$.

We are ready to state and prove the main result of the paper, which is a nonlinear
Bismut-Elworthy formula as the one in theorem \ref{teoBismutLip}, but in the case
of quadratic generator.
\begin{theorem}\label{teoBismut}
 Assume that hypotheses \ref{ip_forward}, \ref{ip_forward_agg}, 
\ref{ip-psiphi} and \ref{ip_G_invertible} hold true, and assume that $\psi$ is differentiable with
respect to $x$, $y$ and $z$. Let $(X^{t,x},Y^{t,x},Z^{t,x})$ be
the solution of the forward-backward system (\ref{fbsde}) and let $U^{h,t,x}$
be defined in (\ref{def-U}).
Then for $0\leq t\leq s\leq T$, $x,h\in H$
\begin{equation}\label{Bismut}
 \E\left[ \nabla_x\,Y^{t,x}_sh \right]=
\E\int_s^T\psi\left(r,X_r^{t,x},Y_r^{t,x},Z_r^{t,x}\right)U^{h,t,x}_r\,dr
+\E\left[  \phi(X_T^{t,x})U^{h,t,x}_T\right].
\end{equation}
\end{theorem}
\dim
We start by approximating the final datum $\phi$
with its inf-sup convolution $\phi_n$ defined in
(\ref{infsupconvphi}).  For all $n\geq 1$ we denote by
 $(Y^{n,t,x},Z^{n,t,x})$ the solution of the Markovian BSDE in
(\ref{fbsde}) with final datum $\phi_n$ 
in the place of $\phi$:
\begin{equation}\label{bsde-n}
    \left\{\begin{array}{l}\dis
 dY^{n,t,x}_\tau=-\psi(\tau,X_\tau,Y^{n,t,x}_\tau,Z^{n,t,x}_\tau)\;d\tau+Z^{n,t,x}_\tau\;dW_\tau,
  \\\dis
  Y^{n,t,x}_T=\phi_n(X_T).
\end{array}\right.
\end{equation}
By theorem \ref{teo_fey_kac}, estimate (\ref{stimaZ-diffle}),
we get that for any $n\geq 1$ there exists a constant $C(n)$, which depends on $n$, such that $C(n)$
is bounded for every $n$ and blows up as $n\rightarrow\infty$, and it is such that
\begin{equation}\label{stimaZ-difflen}
\vert Z_s^{t,x}\vert \leq C(n).
\end{equation}
We get
that the generator $\psi$ acts as a lipschitz generator with respect to $z$
in the BSDE (\ref{bsde-n}), indeed for every $z_1,\,z_2\,\in\Xi$ with $\vert z_i\vert \leq C(n),\,i=1,2$
\[
\vert \psi(s,x,y,z_1)- \psi(s,x,y,z_2)\vert\leq C(n)\vert z_1-z_2\vert.
\]
So for the BSDE (\ref{bsde-n}) the Bismut-Elworthy formula stated in theorem \ref{teoBismutLip}
holds true
\begin{equation}\label{Bismut-n}
 \E\left[ \nabla_x\,Y^{n,t,x}_sh \right]=
\E\int_s^T\psi\left(r,X_r^{t,x},Y_r^{n,t,x},Z_r^{n,t,x}\right)U^{h,t,x}_r\,dr
+\E\left[ \phi_n(X_T^{t,x})U^{h,x}_T\right].
\end{equation}
We have to take the limit as $n\rightarrow\infty.$ We start by computing the limit of the right hand side.
It is immediate to see that
\begin{equation*}
 \lim_{n\rightarrow\infty}\E\left[  \phi_n(X_T^{t,x})U^{h,x}_T\right]=
\E\left[  \phi(X_T^{t,x})U^{h,x}_T\right].
\end{equation*}
Indeed 
\begin{align*}
  \E\vert \left[ \phi_n(X_T^{t,x})-\phi(X_T^{t,x})\right]U^{h,x}_T\vert&\leq 
\left(\E\vert \phi_n(X_T^{t,x})-\phi(X_T^{t,x})\vert^2\right)^{1/2}
\left(\E\vert U^{h,x}_T\vert^2\right)^{1/2}\\
&\leq C \left( T-t\right)^{-1/2}\Vert \phi_n-\phi\Vert_\infty,
\end{align*}
where in the last passage, besides property of the inf-sup convolution,
we have used (\ref{bound-U}) in order to estimate the process $U^{h,x}$.
Now we have to compute
\[
 \lim_{n\rightarrow\infty}\E\int_s^T\psi\left(r,X_r^{t,x},Y_r^{n,t,x},Z_r^{n,t,x}\right)U^{h,t,x}_r\,dr.
\]
To this aim we will show that
\[
 \lim_{n\rightarrow\infty}\E\int_t^T\vert\psi\left(r,X_r^{t,x},Y_r^{n,t,x},Z_r^{n,t,x}\right)U^{h,t,x}_r
-\psi\left(r,X_r^{t,x},Y_r^{t,x},Z_r^{t,x}\right)U^{h,t,x}_r\vert\,dr=0
\]
so that we will deduce that for every $s\in[t,T]$
\[
 \E\int_s^T\psi\left(r,X_r^{t,x},Y_r^{n,t,x},Z_r^{n,t,x}\right)U^{h,t,x}_r\,dr\rightarrow
\E\int_s^T\psi\left(r,X_r^{t,x},Y_r^{t,x},Z_r^{t,x}\right)U^{h,t,x}_r\,dr.
\]
We split the integral with respect to time into two integrals:
\begin{align*}
 &\E\int_t^T\vert\psi\left(r,X_r^{t,x},Y_r^{n,t,x},Z_r^{n,t,x}\right)U^{h,t,x}_r
-\psi\left(r,X_r^{t,x},Y_r^{t,x},Z_r^{t,x}\right)U^{h,t,x}_r\vert\,dr\\ \nonumber
&=\E\int_t^{\frac{t+T}{2}}\vert\psi\left(r,X_r^{t,x},Y_r^{n,t,x},Z_r^{n,t,x}\right)U^{h,t,x}_r
-\psi\left(r,X_r^{t,x},Y_r^{t,x},Z_r^{t,x}\right)U^{h,t,x}_r\vert\,dr\\ \nonumber
&+\E\int_{\frac{t+T}{2}}^T\vert\psi\left(r,X_r^{t,x},Y_r^{n,t,x},Z_r^{n,t,x}\right)U^{h,t,x}_r
-\psi\left(r,X_r^{t,x},Y_r^{t,x},Z_r^{t,x}\right)U^{h,t,x}_r\vert\,dr=I+II.
\end{align*}
We start by estimating I: we recall that by proposition \ref{prop-aprioriZ}, estimate (\ref{stimabismut}), and since
$\Vert \phi_n\Vert_\infty\leq \Vert \phi\Vert_\infty$, there exists a constant $C$, not depending on $n$, such that
\begin{equation}\label{stimabismutn}
\vert Z^{n,t,x}_t\vert\leq C (T-t)^{-1/2}.
\end{equation}
So
\begin{align*}
 I&=\E\int_t^{\frac{t+T}{2}}\vert\psi\left(r,X_r^{t,x},Y_r^{n,t,x},Z_r^{n,t,x}\right)U^{h,t,x}_r
-\psi\left(r,X_r^{t,x},Y_r^{t,x},Z_r^{t,x}\right)U^{h,t,x}_r\vert\,dr\\ \nonumber
 & \leq  \E\int_t^{\frac{t+T}{2}}\left(\vert Z_r^{n,t,x}-Z_r^{t,x}\vert\left(1+\vert Z_r^{n,t,x}\vert
+\vert Z_r^{t,x}\vert \right)\vert U^{h,t,x}_r\vert+\vert Y_r^{n,t,x}-Y_r^{t,x}\vert \vert U^{h,t,x}_r\vert\right)\,dr\\ \nonumber
  &\leq \E\sup_{s\in[t,\frac{t+T}{2}]}\vert Y_s^{n,t,x}-Y_s^{t,x}\vert\int_t^{\frac{t+T}{2}}\vert U^{h,t,x}_r\vert\,dr\\ \nonumber
 +&\sup_{s\in[t,\frac{t+T}{2}]}\vert Z_s^{n,t,x}-Z_s^{t,x}\vert^{1/2}
  \left(1+\vert Z_s^{n,t,x}\vert+\vert Z_s^{t,x}\vert \right)\E\int_t^{\frac{t+T}{2}}\vert Z_r^{n,t,x}-Z_r^{t,x}\vert^{1/2}\vert U^{h,t,x}_r\vert\,dr\\ \nonumber
  &\leq \left(\E\sup_{s\in[t,\frac{t+T}{2}]}\vert Y_s^{n,t,x}-Y_s^{t,x}\vert^3\right)^{1/3}\left(\E\int_t^{\frac{t+T}{2}}\vert U^{h,t,x}_r\vert^{4/3}\,dr\right)^{3/4}\\ \nonumber
 &+C \left(\dfrac{T-t}{2}\right)^{-1/4}\left(\dfrac{T-t}{2}\right)^{-1/2}
  \left(\E\int_t^{\frac{t+T}{2}}\vert Z_r^{n,t,x}-Z_r^{t,x}\vert^{2}\,dr\right)^{1/4}
  \left(\E\int_t^{\frac{t+T}{2}}\vert U^{h,t,x}_r\vert^{4/3}\,dr\right)^{3/4}.
\end{align*}
We estimate the last integral
\begin{align*}
& \left(\E\int_t^{\frac{t+T}{2}}\vert U^{h,t,x}_r\vert^{4/3}\,dr\right)^{3/4}\\ \nonumber
& =\left(\E\int_t^{\frac{t+T}{2}}\vert\dfrac{1}{r-t}\int_t^r \<G^{-1}(s,X^{t,x}_s)\nabla_xX_s^{t,x}h,dW_s \>\vert^{4/3}\,dr\right)^{3/4}\\ \nonumber
 & =\left(\int_t^{\frac{t+T}{2}}\E\vert\dfrac{1}{r-t}\int_t^r \<G^{-1}(s,X^{t,x}_s)\nabla_xX_s^{t,x}h,dW_s \>\vert^{4/3}\,dr\right)^{3/4}\\ \nonumber
 & \leq C\left(\int_t^{\frac{t+T}{2}}\dfrac{1}{(r-t)^{4/3}}\E\left(\int_t^r \vert G^{-1}(s,X^{t,x}_s)\nabla_xX_s^{t,x}h\vert^2\,dr \>\right)^{2/3}\,dr\right)^{3/4}\\ \nonumber
 & \leq C\left(\int_t^{\frac{t+T}{2}}\dfrac{1}{(r-t)^{4/3}}
 \left(r-t \right)^{2/3}\,dr\right)^{3/4}=C\left( T-t \right)^{1/4},\\ \nonumber
\end{align*}
where in the last passage we have used estimate (\ref{stima-nablaX}) for the boundedness of $\nabla_xX^{t,x}$.
Putting together all these estimates we get
\begin{align*}
 I&=\E\int_t^{\frac{t+T}{2}}\vert\psi\left(r,X_r^{t,x},Y_r^{n,t,x},Z_r^{n,t,x}\right)U^{h,t,x}_r
-\psi\left(r,X_r^{t,x},Y_r^{t,x},Z_r^{t,x}\right)U^{h,t,x}_r\vert\,dr\\ \nonumber
&\leq C \left( T-t \right)^{1/4}\left(\E\sup_{s\in[t,\frac{t+T}{2}]}\vert Y_s^{n,t,x}-Y_s^{t,x}\vert^3\right)^{1/3} \\ \nonumber
 &+C \left(\dfrac{T-t}{2}\right)^{-1/4}\left(\dfrac{T-t}{2}\right)^{-1/2}\left( T-t \right)^{1/4}\left(\E\int_t^{\frac{t+T}{2}}
 \vert Z_r^{n,t,x}-Z_r^{t,x}\vert^{2}\,dr\right)^{1/4}\\ \nonumber
 &\leq C \left(\dfrac{T-t}{2}\right)^{-1/2}
  \left(\left(\E\sup_{s\in[t,\frac{t+T}{2}]}\vert Y_s^{n,t,x}-Y_s^{t,x}\vert^3\right)^{1/3}+\left(\E\int_t^{\frac{t+T}{2}}\vert Z_r^{n,t,x}-Z_r^{t,x}\vert^{2}\,dr\right)^{1/4}\right)\rightarrow 0 
\end{align*}
as  $n\rightarrow \infty$, since by well known results in the literature of quadratic BSDEs, as well as a special case of
proposition \ref{prop-convp-bsdequadr}, $Y^{n,t,x}\rightarrow Y^{t,x}$ in $\cals^2$ and it is bounded and $Z^{n,t,x}\rightarrow Z^{t,x}$ in $\calm^2$.

\noindent Next we estimate II: we will use (\ref{bound-Uq-sup})
\begin{align*}
 II&=\E\int_{\frac{t+T}{2}}^T\vert\psi\left(r,X_r^{t,x},Y_r^{n,t,x},Z_r^{n,t,x}\right)U^{h,t,x}_r
-\psi\left(r,X_r^{t,x},Y_r^{t,x},Z_r^{t,x}\right)U^{h,t,x}_r\vert\,dr\\
&\leq\E\sup_{s\in[\frac{t+T}{2},T]}\vert U^{h,t,x}_s\vert \int_{\frac{t+T}{2}}^T\vert\psi\left(r,X_r^{t,x},Y_r^{n,t,x},Z_r^{n,t,x}\right)
-\psi\left(r,X_r^{t,x},Y_r^{t,x},Z_r^{t,x}\right)\vert\,dr\\
 &\leq\left(\E\sup_{s\in[\frac{t+T}{2},T]}\vert U^{h,t,x}_s\vert^q\right)^{1/q}\\
&
\left[\E\left(\int_{\frac{t+T}{2}}^T \vert\psi\left(r,X_r^{t,x},Y_r^{n,t,x},Z_r^{n,t,x}\right)
 -\psi\left(r,X_r^{t,x},Y_r^{t,x},Z_r^{t,x}\right)\vert\,dr\right)^p\right]^{1/p}\\
 &\leq C\dfrac{1}{(T-t)^{1/2}}\left(\E\left(\int_{\frac{t+T}{2}}^T\left(
 \vert Y_r^{n,t,x}-Y_r^{t,x}\vert+
 \vert Z_r^{n,t,x}-Z_r^{t,x}\vert\left(1+\vert Z_r^{n,t,x}\vert+\vert Z_r^{t,x}\vert \right)\right)\,dr\right)^p\right)^{1/p}\\
 &\leq C\dfrac{1}{(T-t)^{1/2}}\left[\dfrac{T-t}{2}
\E\sup_{r\in[t,T]}\vert Y_r^{n,t,x}-Y_r^{t,x}\vert  \right.\\ 
&\left.+\left(\E\left(\int_{\frac{t+T}{2}}^T
 \vert Z_r^{n,t,x}-Z_r^{t,x}\vert^2\,dr\right)^{p/2}
 \left(\int_{\frac{t+T}{2}}^T\left(1+\vert Z_r^{n,t,x}\vert+\vert Z_r^{t,x}\vert \right)^2\,dr\right)^{p/2}\right)^{1/p}\right]\\
 &\leq C\dfrac{1}{(T-t)^{1/2}}\left[\dfrac{T-t}{2}
\E\sup_{r\in[t,T]}\vert Y_r^{n,t,x}-Y_r^{t,x}\vert  \right.\\
&+\left(\E\left(\int_{\frac{t+T}{2}}^T
 \vert Z_r^{n,t,x}-Z_r^{t,x}\vert^2\,dr\right)^{p}\right)^{\frac{1}{2p}}\left(\E
 \left(\int_{\frac{t+T}{2}}^T\left(1+\vert Z_r^{n,t,x}\vert+\vert Z_r^{t,x}\vert \right)^2\,dr\right)^{p}\right)^{\frac{1}{2p}}\rightarrow 0\\
\end{align*}
as $n\rightarrow \infty$. Indeed, by theorem \ref{teo-stimap-bsdequadr},
$Z^{n,t,x}$ as well $Z^{t,x}$ is bounded  in $\calm^{2p}$ by a constant independent
on $n$, and moreover by proposition \ref{prop-convp-bsdequadr} $Z^{n,t,x}$
converges to $Z^{t,x}$ in $\calm^{2p}$.
So we have shown the convergence of $I$ and $II$, and also of the term related to $\phi$
from which we deduce that for every $s\in[t,T]$
\begin{equation}\label{eq}
\lim_{n\rightarrow \infty} \E\left[ \nabla_x\,Y^{n,t,x}_sh \right]=
\E\int_s^T\psi\left(r,X_r^{t,x},Y_r^{t,x},Z_r^{t,x}\right)U^{h,t,x}_r\,dr
+\E\left[ \phi(X_T^{t,x})U^{h,t,x}_T\right].
\end{equation}
In particular, by taking $s=t$ in (\ref{eq}),
\begin{equation*}
\lim_{n\rightarrow \infty}  \nabla_x\,Y^{n,t,x}_th =
\E\int_t^T\psi\left(r,X_r^{t,x},Y_r^{t,x},Z_r^{t,x}\right)U^{h,t,x}_r\,dr
+\E\left[ \phi(X_T^{t,x})U^{h,t,x}_T\right].
\end{equation*}
so we deduce that $\lim_{n\rightarrow \infty}  \nabla_x\,Y^{n,t,x}_th $ exists. Moreover we notice that
this limit is linear in $h$, and we denote it
by $F(t,x)h$,
and moreover we also notice that for every $t\in[0,T]$ and $h\in H$, the map $x\mapsto F(t,x)h$
is continuous.
It remains to show that
\[
 \lim_{n\rightarrow \infty} \nabla_x\,Y^{n,t,x}_th 
=  \nabla_x\,Y^{t,x}_th .
\]
Indeed, for every $\varepsilon>0$ and $x\in H$,
\[
 \dfrac{Y^{n,t,x+\varepsilon h}_t-Y^{n,t,x}_t}{\varepsilon}=\int_0^1Y^{n,t,x+\lambda\varepsilon h}_t h\,d\lambda.
\]
Since
\[
 \dfrac{Y^{n,t,x+\varepsilon h}_t-Y^{n,t,x}_t}{\varepsilon}
\rightarrow \dfrac{Y^{t,x+\varepsilon h}_t-Y^{t,x}_t}{\varepsilon}
\]
and
\[
 \int_0^1Y^{n,t,x+\lambda\varepsilon h}_t h\,d\lambda
\rightarrow \int_0^1 F\left(t,x+\lambda\varepsilon h\right) h\,d\lambda,
\]
and by letting $\varepsilon\rightarrow 0$, we get
\[
F(t,x)=\nabla_x\,Y^{t,x}_t
\]
for almost all $t\in[0,T)$.
Since
\begin{align*}
  \nabla_x\,Y^{t,x}_th =\E\int_t^T\psi\left(r,X_r^{t,x},Y_r^{t,x},Z_r^{t,x}\right)U^{h,t,x}_r\,dr
+\E\left[ \phi(X_T^{t,x})U^{h,t,x}_T\right].
\end{align*}

\noindent The proof is concluded, by noticing that
by the identification \ref{identif-Z}, we also arrive at
\[
 \lim_{n\rightarrow \infty} \E\left[ \nabla_x\,Y^{n,t,x}_sh \right]
= \E\left[ \nabla_x\,Y^{t,x}_sh \right].
\]
for all $0\leq t\leq s<T$.
\qed

We state two corollaries: the first one is about estimates on $\nabla_x\,Y^{t,x}$, and the second one
is about the identification of $\nabla_x\,Y^{t,x}$ with $Z^{t,x}$.
\begin{corollary}\label{cor-stima-nablaY}
 Under the assumptions of theorem \ref{teoBismut},
there exists a constant $C$ depending only on $L_\psi,\,K_\psi,\,K_\phi$ and
on the coefficients of the forward equation (\ref{forward}) such that
\begin{equation}\label{stima-nablaY}
 \vert \nabla_x\,Y^{t,x}\vert \leq C (T-t)^{-1/2}
\end{equation}
\end{corollary}
\dim We start from the Bismut-Elworthy formula we have proved in theorem \ref{teoBismut},
formula (\ref{Bismut}). We first notice that, by the Cauchy-Schwartz inequality and by estimate (\ref{bound-U}),
\begin{equation}\label{stima-nablaY1}
 \E\left[ \phi(X_T^{t,x})U^{h,t,x}_T\right]\leq \Vert\phi\Vert_\infty (T-t)^{-1/2}=K_\phi (T-t)^{-1/2}.
\end{equation}
Next we estimate
\begin{align*}
&\vert \E\int_s^T\psi\left(r,X_r^{t,x},Y_r^{t,x},Z_r^{t,x}\right)U^{h,t,x}_r\,dr\vert\\
&\leq \E\int_t^{\frac{t+T}{2}}\vert\psi\left(r,X_r^{t,x},Y_r^{t,x},Z_r^{t,x}\right)U^{h,t,x}_r
\vert\,dr\\
&+\E\int_{\frac{t+T}{2}}^T\vert\psi\left(r,X_r^{t,x},Y_r^{t,x},Z_r^{t,x}\right)U^{h,t,x}_r
\vert\,dr=I+II\\
\end{align*}
We start by estimating I: by proposition \ref{prop-aprioriZ}, estimate (\ref{stimabismut}), we get
\begin{align}\label{stimaIcor}
 I& \leq C \E\int_t^{\frac{t+T}{2}}\left(1+\vert Y_r^{t,x}\vert+\vert Z_r^{t,x}\vert^2 \right)
\vert U^{h,t,x}_r\vert\,dr\\ \nonumber
 &\leq C \sup_{s\in[t,\frac{t+T}{2}]}
 \left(1+\vert Y_r^{t,x}\vert+\vert Z_s^{t,x}\vert^{3/2} \right)\E\int_t^{\frac{t+T}{2}}
\left(1+\vert Z_r^{t,x}\vert\right)^{1/2}\vert U^{h,t,x}_r\vert\,dr\\ \nonumber
 &\leq C \left(1+\dfrac{T-t}{2}\right)^{-3/4} \left(\E\int_t^{\frac{t+T}{2}}
\left(1+\vert Z_r^{t,x} \vert\right)^{2}\,dr\right)^{1/4}
 \left(\E\int_t^{\frac{t+T}{2}}\vert U^{h,t,x}_r\vert^{4/3}\,dr\right)^{3/4}\\ \nonumber
& \leq C \left(1+\dfrac{T-t}{2}\right)^{-3/4}
\left(\E\int_t^{\frac{t+T}{2}}\vert\dfrac{1}{r-t}\int_t^r \<G^{-1}(s,X^{t,x}_s)\nabla_xX_s^{t,x}h,dW_s \>\vert^{4/3}\,dr\right)^{3/4}\\ \nonumber
 & \leq C\left(1+\dfrac{T-t}{2}\right)^{-3/4}
\left(\int_t^{\frac{t+T}{2}}\dfrac{1}{(r-t)^{4/3}}\E\left(\int_t^r \vert G^{-1}(s,X^{t,x}_s)\nabla_xX_s^{t,x}h\vert^2\,dr \>\right)^{2/3}\,dr\right)^{3/4}\\ \nonumber
 & \leq C\left(1+\dfrac{T-t}{2}\right)^{-3/4}\left( T-t \right)^{1/4}=C\left(\dfrac{T-t}{2}\right)^{-1/2},
\end{align}
where $C$ is a constant depending on $A,\,F,\,G$
 and on $K_\psi,\,K_\phi$ and on $L_\psi$ which gives the linear growth with respect to $y$
and the quadratic growth with respect to $z$ of the generator $\psi$.

\noindent Next we estimate II, and among others we use Lemma \ref{lemma:boundU-q}, estimate (\ref{bound-Uq}):
\begin{align*}
 II=&\E\int_{\frac{t+T}{2}}^T\vert\psi\left(r,X_r^{t,x},Y_r^{t,x},Z_r^{t,x}\right)U^{h,t,x}_r\vert\,dr\\
&\leq\E\sup_{s\in[\frac{t+T}{2},T]}\vert U^{h,t,x}_s\vert \int_{\frac{t+T}{2}}^T\vert\psi\left(r,X_r^{t,x},Y_r^{t,x},Z_r^{t,x}\right)\vert\,dr\\
 &\leq\left(\E\sup_{s\in[\frac{t+T}{2},T]}\vert U^{h,t,x}_s\vert^q\right)^{1/q}\left(\E\left(\int_{\frac{t+T}{2}}^T
 \vert\psi\left(r,X_r^{t,x},Y_r^{t,x},Z_r^{t,x}\right)\vert\,dr\right)^p\right)^{1/p}\\
 &\leq C\dfrac{1}{(T-t)^{1/2}}\left(\E\left(\int_{\frac{t+T}{2}}^T
\left(1+\vert Y_r^{t,x}\vert+\vert Z_r^{t,x}\vert^2 \right)\,dr\right)^p\right)^{1/p}\\
 &\leq C\left(1+\dfrac{1}{(T-t)^{1/2}}\right)
\end{align*}
since by theorem \ref{teo-stimap-bsdequadr} $Y_r^{t,x}$ is bounded and
$Z^{t,x}$ is bounded  in $\calm^{2p}$ by a constant
depending on $K_\psi,\,K_\phi$ and on $L_\psi$ which gives the linear growth with respect to $y$
and the quadratic growth with respect to $z$ of the generator $\psi$.
By putting together the estimate on $II$ with estimate (\ref{stima-nablaY1}) we arrive at
\ref{stima-nablaY} and the proof of the corollary is concluded.
\qed

We notice that differentiability assumptions of $\psi$ with respect to its arguments are needed to achieve the Bismut formula, but do not apper in \ref{stima-nablaY}: this fact will be crucial in the next section, where by means of estimate (\ref{Bismut}) we will solve a semilinear Kolmogorov equation, removing differentiability assumptions on $\psi$.
\begin{corollary}\label{cor-identif-Z}
  Under the assumptions of theorem \ref{teoBismut}, $\forall t\in[0,T]$, $x\in H$
\begin{equation}\label{identif-Z}
  Z^{t,x}_t=\nabla_x\,Y^{t,x}_t G(t).
\end{equation}
\end{corollary}
\dim Let $\phi$ be approximated by its inf-sup convolutions $\phi_n$, and let $(Y^{n,t,x},Z^{n,t,x})$
be the solution of the BSDE (\ref{bsde-n}) with final datum $\phi_n$.
By theorem \ref{teo_fey_kac}, we already know that 
$Z_t^{n,t,x} =\nabla_x Y^{n,t,x}_t G(t)$.
We have just shown in theorem \ref{teoBismut} and corollary \ref{cor-stima-nablaY} that $x\rightarrow Y_\tau^{t,x}$ is differentiable and that $\nabla_xY^{n,t,x}_\tau\rightarrow \nabla_x Y^{t,x}_\tau$,
 $dt\times d\P$ a.e. and a.s.. Moreover, by computing the joint quadratic variation between the process
$v^n(\tau,X_\tau^{t,x}):=Y^{n,t,x}_\tau,\,t\leq\tau\leq T$, and $\int_t^\cdot \xi_s\,dW_s,\,\xi\in L^2_\calp(\Omega\times [0,T])$,
it turns out that
\[
 \int_t^\tau \nabla v^n(s,X_s^{t,x})G(s)\xi_s\,ds
=\int_t^\tau Z^{n,t,x}_s\xi_s\,ds, \;\P\text{ a.s. and for almost all }0\leq \tau\leq T.
\]
By taking a subsequence (that for simplicity we call again $n$) and
letting $n\rightarrow\infty $ in both sides, we get,
\[
 \int_t^\tau \nabla v(s,X_s^{t,x})G(s)\xi_s\,ds
=\int_t^\tau Z^{t,x}_s\xi_s\,ds, \;\P\text{ a.s. and for almost all }0\leq \tau\leq T.
\]
which gives the desired identification,
and as a consequence
\[
 Z^{t,x}_t=  \nabla_x\,Y^{t,x}_t G(t).
\]

\qed

\section{The Bismut formula and mild solutions of a semilinear Kolmogorov equation in the quadratic case}
\label{sez-BismutPDE}
In this section we apply the Bismut formula obtained in theorem \ref{teoBismut}
to solve the semilinear Kolmogorov equation in $H$ given by (\ref{Kolmo}), with $\phi$
and $\psi $ not necessarily differentiable.
We state and prove the main result of this section about the existence and uniqueness of a mild solution to equation (\ref{Kolmo}).
\begin{theorem}\label{teoKolmo}
 Assume that hypotheses \ref{ip_forward}, \ref{ip_forward_agg}, \ref{ip-psiphi} and \ref{ip_G_invertible}
 hold true and let $(X^{t,x},Y^{t,x},Z^{t,x})$ be
the solution of the forward-backward system (\ref{fbsde}) and let $U^{h,t,x}$
be defined in (\ref{def-U}).
Then there exists a unique mild solution $v(t,x)$ of the semilinear Kolmogorov equation (\ref{Kolmo})
given by the formula
\[
 v(t,x)=Y^{t,x}_t,
\]
and such that
\begin{equation}\label{stime-teoKolmo}
\vert v(t,x)\vert\leq C,\qquad \vert \nabla_x v(t,x)\vert\leq C\left(T-t\right)^{-1/2},
\end{equation}
with
\begin{equation}\label{Bismut-teoKolmo}
 \E\left[ \nabla_x\,v(t,x)h\right]=
\E\int_t^T\psi\left(r,X_r^{t,x}, v(r,X_r^{t,x}),\nabla_x v(r,X_r^{t,x})G(r)\right)U^{h,t,x}_r\,dr
+\E\left[  \phi(X_T^{t,x})U^{h,t,x}_T\right].
\end{equation}
\end{theorem}
\dim $Existence$. We start from the case of $\phi $ bounded and continuous as required in the assumptions
of the present theorem, and of $\psi$ also differentiable. We let $\phi_n$ be the inf-sup convolution
of $\phi$, as introduced in (\ref{infsupconvphi}), and we let $(Y^{n,t,x},Z^{n,t,x})$ be the solution
of the BSDE (\ref{bsde-n}) with final datum $\phi_n$, and $v^n$ be the mild solution of a Kolmogorov equation
like (\ref{Kolmo}) with final datum $\phi_n$ instead of $\phi$. Namely $v^n$ satisfies
\begin{equation}
v^n(t,x)=P_{t,T}\left[  \phi_n\right]  \left(  x\right)  +\int_{t}^{T}%
P_{t,s}\left[  \psi(s,\cdot, v^n(s,\cdot), \nabla v^n\left(  s,\cdot\right)G(s)) 
\right]  \left(  x\right)  ds.
\text{\ \ }t\in\left[  0,T\right]  ,\text{ }x\in H. \label{kolmo-mild-n}%
\end{equation}
Since $\phi_n$ is differentiable, by theorem \ref{teo_fey_kac} we already know that
$$
v^n(t,x)=Y^{n,t,x}_t, \qquad\nabla_x v^n(t,x) G(t)=Z^{n,t,x}_t,
$$
moreover
\[
 \lim_{n\rightarrow\infty} Y^{n,t,x}_t=Y^{t,x}_t
\]
where $(Y^{t,x},Z^{t,x})$ is a solution to the backward equation in (\ref{fbsde}). By theorem
\ref{teoBismut} we know that $Y$ satisfies (\ref{Bismut-teoKolmo}) and by corollaries \ref{cor-stima-nablaY}
and \ref{cor-identif-Z} we get respectively estimate \ref{stime-teoKolmo} and the identification
of $\nabla_x v(t,x)G(t)$ with $Z^{t,x}$.
We have to remove differentiability assumptions on $\psi$. To this aim we approximate $\psi$
with its inf-sup convolution $\psi_k$ with respect to $x,\,y$ and $z$, as given in (\ref{infsupconvpsi}), 
We consider the solution of the Kolmogorov equation with nonlinear term given by $\psi_k$ instead of $\psi$,
namely
\begin{equation}
v^k(t,x)=P_{t,T}\left[  \phi\right]  \left(  x\right)  +\int_{t}^{T}%
P_{t,s}\left[  \psi_k(s,\cdot, v^k(s,\cdot), \nabla v^k\left(  s,\cdot\right)G(s)) 
\right]  \left(  x\right)  ds.
\text{\ \ }t\in\left[  0,T\right]  ,\text{ }x\in H. \label{kolmo-mild-k}%
\end{equation}
By the previous part we know that the semilinear Kolmogorov equation (\ref{kolmo-mild-k})
admits a mild solution identified with $Y^{k,t,x}_t$, where $Y^{k,t,x}$ is solution of a backward equation like the one in
the forward-backward system with generator $\psi_k$ instead of $\psi$, namely
\begin{equation}\label{bsde-k}
    \left\{\begin{array}{l}\dis
 dY^{k,t,x}_\tau=-\psi_k(\tau,X_\tau,Y^{k,t,x}_\tau,Z^{k,t,x}_\tau)\;d\tau+Z^{k,t,x}_\tau\;dW_\tau,
  \\\dis
  Y^{k,t,x}_T=\phi(X_T),
\end{array}\right.
\end{equation}
It is well known that $(Y^{k,t,x},Z^{k,t,x})$ converges to $(Y^{t,x},Z^{t,x})$ in
$\cals^2\times\calm^2$, and by proposition \ref{prop-psiconvp-bsdequadr} we know that
$(Y^{k,t,x},Z^{k,t,x})$ converges to $(Y^{t,x},Z^{t,x})$ also in
$\cals^p\times\calm^p$, for any $p\geq 2$. We deduce that
$v^k(t,x)$ converges to $v(t,x)$ and that $v^k(\tau, X_\tau^{t,x})$ converges to $Z_\tau^{t,x}$
in $\calm^p$, and, by taking a subsequence, $dt\times d\P$ a.e..
Next we have to show that the representation (\ref{Bismut-teoKolmo})
holds true also removing differentiability assumptions on $\psi$.
For any $k\geq 1$ it holds true
\begin{equation}\label{Bismut-k}
 \E\left[ \nabla_x\,Y^{k,t,x}_sh \right]=
\E\int_s^T\psi_k\left(r,X_r^{t,x},Y_r^{k,t,x},Z_r^{k,t,x}\right)U^{h,t,x}_r\,dr
+\E\left[ \phi_k(X_T^{t,x})U^{h,x}_T\right].
\end{equation}
We have to take the limit as $k\rightarrow\infty.$ We start by computing the limit of the right hand side, namely we
have to compute
\[
 \lim_{k\rightarrow\infty}\E\int_s^T\psi_k\left(r,X_r^{t,x},Y_r^{k,t,x},Z_r^{k,t,x}\right)U^{h,t,x}_r\,dr.
\]
To this aim we will show that
\[
 \lim_{k\rightarrow\infty}\E\int_t^T\vert\psi_k\left(r,X_r^{t,x},Y_r^{k,t,x},Z_r^{k,t,x}\right)U^{h,t,x}_r
-\psi\left(r,X_r^{t,x},Y_r^{t,x},Z_r^{t,x}\right)U^{h,t,x}_r\vert\,dr=0.
\]
We start by splitting the integral with respect to time into two integrals:
\begin{align*}
 &\E\int_t^T\vert\psi_k\left(r,X_r^{t,x},Y_r^{k,t,x},Z_r^{k,t,x}\right)U^{h,t,x}_r
-\psi\left(r,X_r^{t,x},Y_r^{t,x},Z_r^{t,x}\right)U^{h,t,x}_r\vert\,dr\\ \nonumber
&=\E\int_t^{\frac{t+T}{2}}\vert\psi_k\left(r,X_r^{t,x},Y_r^{k,t,x},Z_r^{k,t,x}\right)U^{h,t,x}_r
-\psi\left(r,X_r^{t,x},Y_r^{t,x},Z_r^{t,x}\right)U^{h,t,x}_r\vert\,dr\\ \nonumber
&+\E\int_{\frac{t+T}{2}}^T\vert\psi_k\left(r,X_r^{t,x},Y_r^{k,t,x},Z_r^{k,t,x}\right)U^{h,t,x}_r
-\psi\left(r,X_r^{t,x},Y_r^{t,x},Z_r^{t,x}\right)U^{h,t,x}_r\vert\,dr=I+II.
\end{align*}
We start by estimating I: we recall that by proposition \ref{prop-aprioriZ}, estimate (\ref{stimabismut}), and by
properties of the inf-sup convolutions, there exists a constant $C$, not depending on $k$, such that
\begin{equation}\label{stimabismutk}
\vert Z^{k,t,x}_t\vert\leq C (T-t)^{-1/2},
\end{equation}
and moreover by corollary \ref{cor-stima-nablaY}, estimate (\ref{stimaIcor})
\[
 \E\int_t^{\frac{t+T}{2}}\left(1+\vert Y_r^{k,t,x}\vert+\vert Z_r^{k,t,x}\vert^2 \right)
\vert U^{h,t,x}_r\vert\,dr\leq C\left(T-t\right)^{-1/2},
\]
where $C$ does not depend on $k$.
So
\begin{align*}
 I&\leq \E\int_t^{\frac{t+T}{2}}\vert\psi_k\left(r,X_r^{t,x},Y_r^{k,t,x},Z_r^{k,t,x}\right)U^{h,t,x}_r
-\psi\left(r,X_r^{t,x},Y_r^{k,t,x},Z_r^{k,t,x}\right)U^{h,t,x}_r\vert\,dr\\
&+\E\int_t^{\frac{t+T}{2}}\vert\psi\left(r,X_r^{t,x},Y_r^{k,t,x},Z_r^{k,t,x}\right)U^{h,t,x}_r
-\psi\left(r,X_r^{t,x},Y_r^{t,x},Z_r^{t,x}\right)U^{h,t,x}_r\vert\,dr\\
 &\leq C\left( \sup_{x\in H,\,y\in\R,\,z \in H} \vert \dfrac{\psi_k(t,x,y,z)}{1+\vert z\vert^2}-
 \dfrac{\psi(t,x,y,z)}{1+\vert z\vert^2}\vert\right) \int_t^{\frac{t+T}{2}}
\left(1+\vert Y_r^{k,t,x}\vert+\vert Z_r^{k,t,x}\vert^2\right)\vert U^{h,t,x}_r\vert\,dr
 \\ \nonumber
  &+C \E\int_t^{\frac{t+T}{2}}\vert Y_r^{k,t,x}-Y_r^{t,x}\vert\vert U^{h,t,x}_r\vert\,dr
 +C \E\int_t^{\frac{t+T}{2}}\vert Z_r^{k,t,x}-Z_r^{t,x}\vert\left(1+\vert Z_r^{k,t,x}\vert+\vert Z_r^{t,x}\vert \right)\vert U^{h,t,x}_r\vert\,dr\\ \nonumber
   &\leq C (T-t)^{-1/2}\left(\sup_{x\in H,\,y\in\R,\,z \in H} \vert \dfrac{\psi_k(t,x,y,z)}{1+\vert z\vert^2}-
   \dfrac{\psi(t,x,y,z)}{1+\vert z\vert^2}\vert\right)
   \\ \nonumber
 &+C \E\sup_{s\in[t,\frac{t+T}{2}]} \vert U^{h,t,x}_s\vert\int_t^{\frac{t+T}{2}}\vert Y_r^{k,t,x}-Y_r^{t,x}\vert\,dr\\
   &+ C\sup_{s\in[t,\frac{t+T}{2}]}\vert Z_s^{k,t,x}-Z_s^{t,x}\vert^{1/2}
   \left(1+\vert Z_s^{k,t,x}\vert+\vert Z_s^{t,x}\vert \right)\E\int_t^{\frac{t+T}{2}}\vert Z_r^{k,t,x}-Z_r^{t,x}\vert^{1/2}\vert U^{h,t,x}_r\vert\,dr\\ \nonumber
   &\leq C (T-t)^{-1/2}\left(\sup_{x\in H,\,y\in\R,\,z \in H} \vert \dfrac{\psi_k(t,x,y,z)}{1+\vert z\vert^2} -
   \dfrac{\psi(t,x,y,z)}{1+\vert z\vert^2}\vert\right)\\
 &+C \left(\E\sup_{s\in[t,\frac{t+T}{2}]} \vert U^{h,t,x}_s\vert^2\right)^{1/2}\left(
 \E\int_t^{\frac{t+T}{2}}\vert Y_r^{k,t,x}-Y_r^{t,x}\vert^2\,dr\right)^{1/2}\\
   &+ C \left(\dfrac{T-t}{2}\right)^{-1/4}\left(\dfrac{T-t}{2}\right)^{-1/2}
   \left(\E\int_t^{\frac{t+T}{2}}\vert Z_r^{k,t,x}-Z_r^{t,x}\vert^{2}\,dr\right)^{1/4}
   \left(\E\int_t^{\frac{t+T}{2}}\vert U^{h,t,x}_r\vert^{4/3}\,dr\right)^{3/4}\\ \nonumber
  &\leq C \left(\dfrac{T-t}{2}\right)^{-1/2}\left[\left(\sup_{x\in H,\,y\in\R,\,z \in H} \vert \dfrac{\psi_k(t,x,y,z)}{1+\vert z\vert^2}-   \dfrac{\psi(t,x,y,z)}{1+\vert z\vert^2}\vert\right)
 \right.\\ \nonumber
    &\left.+\left(\E\int_t^{\frac{t+T}{2}}\vert Y_r^{k,t,x}-Y_r^{t,x}\vert^{2}\,dr\right)^{1/2}
    +\left(\E\int_t^{\frac{t+T}{2}}\vert Z_r^{k,t,x}-Z_r^{t,x}\vert^{2}\,dr\right)^{1/4}\right]\rightarrow 0
\end{align*}
as $k\rightarrow \infty$. Indeed by theorem \ref{teo-stimap-bsdequadr} $Z^{n,t,x}$ is bounded in $\calm^p$, for any $p\geq 1$, and since
by well known results in the literature of quadratic BSDEs, as well as a special case of
proposition \ref{prop-psiconvp-bsdequadr}, $(Y^{n,t,x},Z^{n,t,x})\rightarrow (Y^{t,x},Z^{t,x})$
in $\calm^2$.

\noindent Next we estimate II:
\begin{align*}
 &II\leq\E\int_{\frac{t+T}{2}}^T\vert\psi_k\left(r,X_r^{t,x},Y_r^{k,t,x},Z_r^{k,t,x}\right)U^{h,t,x}_r
-\psi\left(r,X_r^{t,x},Y_r^{k,t,x},Z_r^{k,t,x}\right)U^{h,t,x}_r\vert\,dr\\
&+\E\int_{\frac{t+T}{2}}^T\vert\psi\left(r,X_r^{t,x},Y_r^{k,t,x},Z_r^{k,t,x}\right)U^{h,t,x}_r
-\psi\left(r,X_r^{t,x},Y_r^{t,x},Z_r^{t,x}\right)U^{h,t,x}_r\vert\,dr\\
 &\leq C \sup_{x\in H,\,y\in\R,\,z \in H} \vert \dfrac{\psi_k(t,x,y,z)}{1+\vert z\vert^2}-
   \dfrac{\psi(t,x,y,z)}{1+\vert z\vert^2}\vert
 \E\int_{\frac{t+T}{2}}^T 
 \left(1+\vert Z_r^{k,t,x}\vert^2\right)\vert U^{h,t,x}_r\vert\,dr\\
  \nonumber
 &+\E\sup_{s\in[\frac{t+T}{2},T]}\vert U^{h,t,x}_s\vert \int_{\frac{t+T}{2}}^T\vert\psi\left(r,X_r^{t,x},Y_r^{k,t,x},Z_r^{k,t,x}\right)
 -\psi\left(r,X_r^{t,x},Y_r^{t,x},Z_r^{t,x}\right)\vert\,dr\\
  &\leq C \sup_{x\in H,\,y\in\R,\,z \in H}  \vert\dfrac{\psi_k(t,x,y,z)}{1+\vert z\vert^2}-
   \dfrac{\psi(t,x,y,z)}{1+\vert z\vert^2}\vert
 \E\sup_{s\in[\frac{t+T}{2},T]}\vert U^{h,t,x}_s\vert\int_{\frac{t+T}{2}}^T 
 \left(1+\vert Z_r^{k,t,x}\vert^2\right)\,dr\\
 &+\left(\E\sup_{s\in[\frac{t+T}{2},T]}\vert U^{h,t,x}_s\vert^q\right)^{1/q}
 \left(\E\left(\int_{\frac{t+T}{2}}^T
  \vert\psi\left(r,X_r^{t,x},Y_r^{k,t,x},Z_r^{k,t,x}\right)
  -\psi\left(r,X_r^{t,x},Y_r^{t,x},Z_r^{t,x}\right)\vert\,dr\right)^p\right)^{1/p}\\
   & \leq C \sup_{x\in H,\,y\in\R,\,z \in H} \vert \dfrac{\psi_k(t,x,y,z)}{1+\vert z\vert^2}-
    \dfrac{\psi(t,x,y,z)}{1+\vert z\vert^2}\vert
  \left(\E\sup_{s\in[\frac{t+T}{2},T]}\vert U^{h,t,x}_s\vert^2\right)^{1/2}\\
&\;\;\left(\E\left(\int_{\frac{t+T}{2}}^T 
  \left(1+\vert Z_r^{k,t,x}\vert^2\right)\,dr\right)^2\right)^{1/2}\\
  &+ C\dfrac{1}{(T-t)^{1/2}}\left(\E\left(\int_{\frac{t+T}{2}}^T\left(
   \vert Y_r^{k,t,x}-Y_r^{t,x}\vert+
   \vert Z_r^{k,t,x}-Z_r^{t,x}\vert\left(1+\vert Z_r^{k,t,x}\vert+\vert Z_r^{t,x}\vert \right)\right)\,dr\right)^p\right)^{1/p}\\
    &\leq C \dfrac{1}{(T-t)^{1/2}}\sup_{x\in H,\,y\in\R,\,z \in H} 
\vert \dfrac{\psi_k(t,x,y,z)}{1+\vert z\vert^2}-
    \dfrac{\psi(t,x,y,z)}{1+\vert z\vert^2}\vert
  \left(\E\left(\int_{\frac{t+T}{2}}^T 
  \left(1+\vert Z_r^{k,t,x}\vert^2\right)\,dr\right)^2\right)^{1/2}\\
  &+ C\dfrac{1}{(T-t)^{1/2}}\left[\E\int_{\frac{t+T}{2}}^T\vert Y_r^{k,t,x}-Y_r^{t,x}\vert^p\,dr+\left(\E\left(\int_{\frac{t+T}{2}}^T
   \vert Z_r^{k,t,x}-Z_r^{t,x}\vert^2\,dr\right)^{p/2}\right.\right.\\
  &\left.\left.\left(\int_{\frac{t+T}{2}}^T\left(1+\vert Z_r^{k,t,x}\vert+\vert Z_r^{t,x}\vert \right)^2\,dr\right)^{p/2}\right)^{1/p}\right]
\rightarrow 0\\
\end{align*}
as $k\rightarrow \infty$. Indeed, by theorem \ref{teo-stimap-bsdequadr},
$Z^{k,t,x}$ as well $Z^{t,x}$ is bounded  in $\calm^{2p}$ by a constant independent
on $k$, and moreover by proposition \ref{prop-convp-bsdequadr} $(Y^{k,t,x},Z^{k,t,x})$
converges to $(Y^{t,x},Z^{t,x})$ in $\calm^{2p}$.
So we have shown the convergence of $I$ and $II$,
from which we deduce that for every $s\in[t,T]$
\begin{equation*}
\lim_{k\rightarrow \infty} \E\left[ \nabla_x\,Y^{k,t,x}_sh \right]=
\E\int_s^T\psi\left(r,X_r^{t,x},Y_r^{t,x},Z_r^{t,x}\right)U^{h,t,x}_r\,dr
+\E\left[ \phi(X_T^{t,x})U^{h,x}_T\right].
\end{equation*}
As before we can show that
\[
 \lim_{k\rightarrow \infty} \E\left[ \nabla_x\,Y^{k,t,x}_sh \right]
= \E\left[ \nabla_x\,Y^{t,x}_sh \right].
\]
The identification
\[
 \nabla_x\,Y^{t,x}_t= Z^{t,x}_t G(t).
\]
can be obtained as in corollary \ref{cor-identif-Z}, and the proof is concluded.
\qed

We conclude this Section summing up what type of Kolmogorov equation we are able to treat.
\begin{remark}\label{remarkKolmo}
We notice that under the invertibility assumptions on $G$ we are able to solve a Kolmogorov equation
(\ref{Kolmo}) for the unknown $v$ with non linear term $\psi$ with quadratic growth with respect to the derivative $\nabla v$, and lipschitz continuous with respect to $v$ and $x$, and with final datum only continuous and bounded. We notice that, due to the boundedness of $v$ given by the estimates on $Y$ in Proposition \ref{teo-stimap-bsdequadr}, linear growth with respect to $v$  may be removed with some technical efforts, that we omit here.

Coming to a comparison with the existing literature, we are able to treat a superquadratic Kolmogorov equation with final datum bounded and continuous: here we ask invertibility assumptions on $G(t)$, in \cite{Mas1} a similar result is achieved with $A$ and $G$ commuting, while in \cite{MR} it is considered
a Kolmogorov equation \ref{Kolmo} a final datum  locally lipschitz continuous and not necessarily bounded but with polynomial growth with respect to $x$, so in \cite{MR} the request of regularity on the final are significantly stronger than in the present paper.
\end{remark}

\section{A quadratic optimal control problem}
\label{sez-appl-contr}
Now we apply the above results to perform the synthesis of the optimal control
for a class of control problems with nonlinear state equation and with related current
cost with quadratic growth with respect to the control $u$ and the final cost bounded and only continuous.

Let $X^u$ be the solution of the controlled state equation
\begin{equation}
\left\{
\begin{array}
[c]{l}%
dX^{u}_\tau  =AX^{u}_\tau d\tau +F(\tau,X_\tau)d\tau+R\left( u_\tau \right)  d\tau+
G(\tau)dW_\tau ,\text{ \ \ \ }\tau\in\left[  t,T\right] \\
X^{u}_t  =x.
\end{array}
\right.  \label{sdecontrolforte}%
\end{equation}
Notice that, due to the invertibility assumptions on $G$, equation (\ref{sdecontrolforte})
can be rewritten as
\begin{equation*}
\left\{
\begin{array}
[c]{l}%
dX^{u}_\tau  =AX^{u}_\tau d\tau +F(\tau,X_\tau)d\tau+G(\tau)\tilde R\left( u_\tau \right)  d\tau+
G(\tau)dW_\tau ,\text{ \ \ \ }\tau\in\left[  t,T\right] \\
X^{u}_t  =x,
\end{array}
\right.
\end{equation*}
where $\tilde R\left( u_\tau \right) =G^{-1}(\tau )R\left( u_\tau \right) $. This means that equation
(\ref{sdecontrolforte}) can be written with the ``special structure'' that allows to study
the optimal control problem related by means of BSDEs.
On $R$ we make the following assumption:
\begin{hypothesis}\label{ip-R}
The space $U$ where the control process takes its values is a general Banach space.
For the map $R:U\rightarrow H$ there exists a constant $c>0$ such that
$\vert R(u)\vert \leq c(1+\vert u\vert_U)$.
\end{hypothesis}
Beside equation (\ref{sdecontrolforte}), we define the cost
\begin{equation}
J\left(  t,x,u\right)  =\mathbb{E}\int_{t}^{T}
g\left(s,X^{u}_s,u_s\right)ds+\mathbb{E}\phi\left(X^{u}_T\right). 
\label{cost}%
\end{equation}
for real functions $g$ on $[0,T]\times H\times U$ and  $\mathbb{\phi}$ on $H$.
The control problem in strong formulation is to
minimize this functional $J$ over all admissible controls $u$.
By admissible control we mean  an $(\calf_t)_t$-predictable
process, taking values in a closed subset $K$ of $U$, such that
\[
 \E\int_0^T \vert u_s \vert ^2 ds < + \infty.
\]
This assumption is natural this time since we assume here that the cost 
has quadratic growth at infinity, as it can be seen in the
following assumptions on the cost $J$.

\begin{hypothesis}
\label{ip costo}
We assume:
\begin{enumerate}
\item $\mathbb{\phi}:H\rightarrow\mathbb{R}$ is bounded and continuous;

\item $g:[0,T]\times H\times U\rightarrow\mathbb{R}$ is measurable and for all $t\in[0,T],\,u\in U$,
$x \mapsto g(t,x,u)$ is bounded and continuous, moreover for all $t\in[0,T],\,x\in H,\,u\in U$,
there exists a constant $c>0$ such that
\begin{equation}
 \label{crescita costo1}
0\leq g(t,x,u)\leq c(1+ \vert u \vert ^2) 
\end{equation}
and there exist $R>0$, $C>0$ such that
\begin{equation}
\label{crescita costo}
g(t,x,u)\geq C \vert u \vert ^2 \qquad \forall u \in K,\, \vert u\vert \geq R.
\end{equation}
\item $l$ is lipschitz continuous with respect to $x$, uniformly with respect to $t\in[0,T]$
and $u\in U$, that is for all $t\in[0,T],\,x_1,x_2\in H,\,u\in U$, for some $C>0$,
\begin{equation*}
 \vert g(t,x_1,u)-g(t,x_2,u) \vert\leq C\vert x_1-x_2 \vert.
\end{equation*}
\end{enumerate}
\end{hypothesis}

We define in a classical way the Hamiltonian function relative to the above
problem:%
\begin{equation}\label{hamilton}
\psi\left(t,x,z\right)  =\inf_{u\in K}\left\{  g\left(t,x,u\right)
+zR(u)\right\}\quad \forall z\in H .
\end{equation}
We prove that the Hamiltonian function just defined satisfies the polynomial growth
conditions and the local lipschitzianity required in hypothesis \ref{ip-psiphi}.
\begin{lemma}
 \label{lemma-hamilton}Assume that hypotheses \ref{ip-R} and \ref{crescita costo1}, point 1 and 2, hold true.
Then the Hamiltonian $\psi:[0,T]\times H\times H\rightarrow \R$
is Borel measurable, there exists a constant $C>0$ such that
\[
-C(1+\vert z\vert^2)\leq \psi(t,x,z)  \leq g(t,x,u)+C \vert z\vert (1+\vert u\vert), \qquad \forall u \in K.
\]
Moreover if the infimum in (\ref{hamilton}) is attained, it is attained
in a ball of radius $C(1+\vert z\vert)$,  that is 
\[
\psi(t,x,z)  =\inf_{u\in K,\vert u\vert \leq C(1+\vert z\vert)}\left\{  g\left(t,x,u\right)
+zR(u)\right\},\quad z \in H,
\]
and
\[
\psi(t,x,z)  <  g\left(t,x,u\right)
+zR(u)\quad \text{if }\vert u\vert > C(1+\vert z\vert).
\]
In particular it follows that $\psi$ is locally lipschitz continuous with respect to $z$, namely
for all $t\in[0,T],\,x\in H,\,z_1,z_2\in H$, for some $C>0$,
\begin{equation}\label{stima2psi}
 \vert \psi(t,x,z_1)-\psi(t,x,z_2) \vert\leq C(1+\vert z_1\vert+\vert z_2 \vert)\vert z_1-z_2 \vert.
\end{equation}
Moreover, if hypothesis \ref{crescita costo1}, point 3, holds true, then
$\psi$ is lipschitz continuous with respect to $x$, namely
for all $t\in[0,T],\,x_1,x_2\in H,\,z\in H$, for some $C>0$,
\begin{equation}\label{stimaxpsi}
 \vert \psi(t,x_1,z)-\psi(t,x_2,z) \vert\leq C\vert x_1-x_2 \vert.
\end{equation}
\end{lemma}
\dim The proof is given in \cite{fuhute}, lemma 3.1, apart from (\ref{stima2psi})
and (\ref{stimaxpsi}), that we briefly discuss here.
For what concerns (\ref{stima2psi}), for every $u\in K$, $u$ in the set where the infimum in the definition of the hamiltonian (\ref{hamilton}) is achieved,
we get for all $t\in[0,T],\,x\in H,\,z_1,z_2\in H$
\[
-C(1+\vert z_1\vert +\vert z_2\vert)\vert z_1-z_2\vert\leq 
g\left(t,x,u\right)+z_1u-g\left(t,x,u\right)-z_2u\leq C(1+\vert z_1\vert +\vert z_2\vert)\vert z_1-z_2\vert
\]
which gives (\ref{stima2psi}).

\noindent For what concerns (\ref{stimaxpsi}), if hypothesis
\ref{crescita costo1}, point 3, holds true, then for all $t\in[0,T],\,x_1,x_2\in H,\,z\in H$
\[
 -C\vert x_1-x_2\vert\leq g\left(t,x_1,u\right)-g\left(t,x_2,u\right)\leq  g\left(t,x_1,u\right)+zR(u)
-\inf_{u\in K}\left\{  g\left(t,x_2
,u\right)
+zR(u)\right\}
\]
and since this inequality is true $\forall u\in U$, we immediately get 
\[
\psi (t,x_1,z)-\psi(t,x_2,z)\geq
 -C\vert x_1-x_2\vert.
\]
Arguing in 
a similar way we arrive at
\[
\psi (t,x_1,z)-\psi(t,x_2,z)\leq
 C\vert x_1-x_2\vert.
\]
and this gives (\ref{stimaxpsi}) and
concludes the proof.
\qed

We define
\begin{equation}\label{defdigammagrande}
\Gamma(t,x,z)=\left\{ u\in U: zR(u)+g(t,x,u)= \psi(t,x,z)\right\};
\end{equation}
if $\Gamma(t,x,z) \neq \emptyset$ for every $t\in[0,T], x\in H, z\in H$, by \cite{AuFr}, see Theorems 8.2.10 and
8.2.11, $\Gamma$ admits a measurable selection, i.e. there exists
a measurable function $\gamma:[0,T]\times H\times H \rightarrow U$ with
$\gamma(t,x,z)\in \Gamma(t,x,z)$ for every $t\in[0,T], x\in H, z\in H$.

In the following theorem we will prove the fundamental relation, by applying theorem
\ref{teoKolmo}.

\begin{theorem}\label{th-rel-font}
 Assume hypotheses \ref{ip_forward}, \ref{ip_forward_agg}, \ref{ip_G_invertible}, 
\ref{ip costo} hold true.
 Let $v$ the solution of the HJB equation (\ref{Kolmo}). For
every $t\in [0,T]$, $x\in H$ and for
 all admissible control $u$ we have $J(t,x,u(\cdot))
 \geq v(t,x)$,
  and the
 equality holds if and only if
$$
u_s\in \Gamma\left( t,x,\nabla 
v(s ,X^{u,t,x}_s)G(s)
\right)
  $$
\end{theorem}

{\bf Proof.} The proof follows from proposition 4.1 in \cite{fuhute},
recalling that by theorem \ref{teoKolmo} equation (\ref{Kolmo}) admits a unique mild solution
$v(t,x)=Y_t^{t,x}$, where $(Y^{t,x},Z^{t,x})$ is solution to the BSDE
in FBSDE (\ref{fbsde}), and that by corollary \ref{cor-identif-Z}
$Z_t^{t,x}=\nabla v(t, x)G(t)$.
\qed

With the assumptions of Theorem
\ref{th-rel-font}, we can define the so called
optimal feedback law 
\begin{equation}\label{leggecontrolloottima}
u(s,x)=\gamma\Big(\nabla
v(s ,X^{u,t,x}_s) G(s)\Big),\qquad
s\in [t,T],\;x\in H,
\end{equation}
and the related closed loop equation in mild form is given by
\begin{equation}\label{cleq}
\overline{X}_s= e^{(s-t)A}x+\int_t^s e^{(s-r)A}F(r,\overline{X}_r)\,dr
+\int_{t}^s e^{(s-r)A}R\left(\gamma(r,\overline{X}_r, \nabla v(r, \overline{X}_r)\right)+
\int_t^s e^{(s-r)A}G(r)\,dW_r.
\end{equation}
If the closed loop equation admits a solution the pair $(\overline{u}=u(s,\overline{X}_s),\overline{X}_s)_{s\in[t,T]}$
is optimal for the control problem.
Due to the lack
of regularity of the feedback law $u$ occurring in
(\ref{cleq}), the existence of a solution of the closed loop
equation is not obvious
This problem can be avoided by formulating the optimal control problem in the weak sense

In the following, by an admissible control system
we mean
$$(\Omega,\mathcal{F},
\left(\mathcal{F}_{t}\right) _{t\geq 0}, \mathbb{P}, W,
u(\cdot),X^u),$$
where $W$ is an $H$-valued Wiener process, $u$ is an admissible control and $X^u$
solves the controlled equation (\ref{sdecontrolforte}).
The control problem in weak formulation is to minimize
the cost functional over all the admissible control systems.
\begin{theorem}\label{teo su controllo debole}
Assume hypotheses \ref{ip_forward}, \ref{ip_forward_agg}, \ref{ip_G_invertible}, 
\ref{ip costo} hold true.
 Let $v$ the solution of the HJB equation (\ref{Kolmo}). 
 For
every $t\in [0,T]$, $x\in H$ and for
 all admissible control systems we have $J(t,x,u(\cdot))
 \geq v(t,x)$,
  and the
 equality holds if and only if
$$
u_s\in \Gamma\left( s,X^{u}_s, \nabla
v(s ,X^{u}_s)G(s)
\right)
  $$
Moreover
assume that the set-valued map $\Gamma$ is not empty and let $\gamma$
be its measurable selection. Then the process defined by
\begin{equation*}
u_{\tau}=\gamma(\tau , X^{u}_\tau,\nabla v(\tau,X^u_\tau )G(\tau))
,\text{ \ }\mathbb{P}\text{-a.s. for a.a. }\tau
\in\left[ t,T\right]
\end{equation*}
is optimal.

\noindent Finally, the closed loop equation \ref{cleq}
admits a weak solution
$(\Omega,\mathcal{F},
\left(\mathcal{F}_{t}\right) _{t\geq 0}, \mathbb{P}, W,X)$
which is unique in law and setting
$$
u_{\tau}=\gamma\left(\tau , X^{u}_\tau,\nabla v(\tau,X^u_\tau )G(\tau)\right),
$$
we obtain an optimal admissible control system $\left(
W,u,X\right) $.
\end{theorem}

\noindent {\bf Proof.} The proof follows from the fundamental relation
stated in theorem \ref{th-rel-font};
the closed loop equation can be solved as in \cite{fuhute}, proposition
5.2.
\qed

\subsection{Optimal control problems for a semilinear heat equation}
\label{sez_contr_heat}

In this section we briefly show how to apply our results to
solve the optimal control problem when the
state equation is a general semilinear heat equation with additive noise.

Namely we consider a bounded domain in $\R^n$ denoted by $\calo$,
$H=L^2(\calo)$.
We consider the following controlled heat equation, for $0\leq t\leq s\leq T$,
\begin{equation}\label{heat equation contr}
 \left\{
  \begin{array}{l}
  \dis
\frac{ \partial y}{\partial s}(s,\xi)= \Delta y(s,\xi)+f\left(s,\xi,y(s,\xi)\right)
+ \sigma\left( s,\xi\right) r\left( \xi, u\left(s,\xi\right)\right)+
\sigma\left( s,\xi\right)\frac{ \partial W}{\partial s}(s,\xi), \;
\xi\in \calo,\\
\dis
y(t,\xi)=x(\xi),\;\xi\in \calo,
\\\dis
 y(s,\xi)=0, \quad \xi\in\partial \calo.
\end{array}
\right.
\end{equation}
where $u_s\in L^2(\calo)$ represents the control.
In the following we denote by 
$\mathcal{A}_{d}$ the set of
admissible controls, that is the real valued predictable processes such that 
\[
 \E \int_0^T \left( \int_{\calo}\vert u_t(\xi) \vert^2 d\xi\right)  dt <+\infty.
\]
and such that $u_t \in K$, where $K$ is a closed subset of $H$, not
necessarily coinciding with $H$, where this time $U=H$.
The process $W(s,\xi)$ is a space time white noise on $[0,T]\times \calo$.

\noindent Our aim is to minimize over all admissible controls the cost functional
\begin{equation}
J\left(  t,x(\xi),u\right)  =\mathbb{E}\int_{t}^{T}\int_{\calo}
\bar l\left(s,\xi,y(s,\xi),u_s(\xi)\right)d\xi ds
+\mathbb{E}\int_{\calo} \bar\phi\left(\xi,y(T,\xi)\right)\,d\xi. 
\label{heat cost}%
\end{equation}
for real functions $\bar\phi$ and $\bar l$.

We make the following assumptions on the controlled stochastic heat equation
(\ref{heat equation contr}) and on the related cost $J$.

\begin{hypothesis}
\label{ip costo heat}The functions $f,$ $\sigma,$ $r,$ $\bar l,$ $\bar\phi$
are all Borel measurable and real valued. Moreover,

\begin{enumerate}
\item $f:[ 0,T] \times\calo \times\mathbb{R}%
\longrightarrow\mathbb{R}$ is continuous; for every $s\in[ 0,T]
$ and every $\xi\in\calo$, we have $f( s,\xi ,\cdot) \in
C^{1}( \mathbb{R})$; and there exists $c_{1}$ continuous on $[
0,1] $ such that
\[
\left\vert f\left( s,\xi,x\right) \right\vert \leq c_{1}\left(
\xi\right) \left( 1+\left\vert x\right\vert \right) \text{, \
\ \ }\left\vert \nabla_{x}f\left( s,\xi,x\right) h\right\vert
\leq c_{1}\left( \xi\right) \left\vert h\right\vert 
\]
for every $s\in[ 0,T] ,$ $\xi\in\calo ,$ $x,h\in\mathbb{R}$.

\item $\sigma:[0,T]\times\calo\rightarrow\R$ is bounded and positive, invertible with a bounded inverse;
\item for every $\xi\in\calo,\,u\in \R$, $\vert r(\xi,u)\vert\leq c(1+\vert u\vert)$;
\item $\bar l:[0,T]\times\calo\times\mathbb{R}
\times\R\rightarrow\mathbb{R}$ is continuous and for every $s\in[0,T],
\, \xi\in\calo, x,x_1,X_2\in\R,\,u\in \R$ we have 
\[
0\leq \bar l(s,\xi,x,u) \leq c(1+\vert u\vert^2)
\]
and there exists $R>0$, $C>0$ such that
\begin{equation*}
\bar l(s,\xi,x,u)\geq C \vert u \vert ^2 \qquad  \vert u\vert \geq R \text{ for all }u\in K;
\end{equation*}
moreover
\[
\vert \bar l(s,\xi,x_1,u)-\bar l(s,\xi,x_2,u)\vert \leq \vert x_1-x_2\vert;
\]
\item $\bar\phi:\calo \times\mathbb{R}\rightarrow\mathbb{R}$ is
continuous and bounded
\item $x_{0}\in L^2(\calo)$;
\end{enumerate}
\end{hypothesis}
Let us define, for $s\in[0,T], x\in H, u\in U$
\begin{equation}\label{heat notazioni}
\begin{array}{ll}
 F(s,x)(\xi):=f\left(s,\xi,x(\xi)\right), & \left( G(s)z\right) \left( \xi\right)
=\sigma\left( s,\xi\right) z\left(
\xi\right), \\
\left( Ru\right) \left( \xi\right) =r\left( \xi,u\left(
\xi\right) \right),&\\
l(s,x,u)=\int_{\calo}
\bar l\left(s,\xi,x(\xi, u(\xi))\right)\,d\xi,& \phi (x)=\int_{\calo} \bar\phi\left(\xi,x(\xi)\right)\,d\xi
\end{array}
\end{equation}
It turns out that if $f,\, r,\,\sigma,\,\bar l$ and $\bar\phi$ satisfy
hypothesis \ref{ip costo heat}, then $F,\, R,\,G,\,l$ and $\phi$
defined in (\ref{heat notazioni}) satisfy hypothesis
\ref{ip costo}.
\noindent Moreover equation (\ref{heat equation contr}) can be written in an abstract way in $H$
as
\begin{equation}\label{heat eq abstr contr}
 \left\{
\begin{array}
[c]{l}%
dX^u_\tau  =AX^u_\tau d\tau+F(\tau, X_\tau^u)d\tau+R(u_\tau) d\tau+G(\tau)dW_\tau
,\text{ \ \ \ }\tau\in\left[  t,T\right] \\
X^u_t =x_0,
\end{array}
\right.
\end{equation}
where $A$ is the Laplace operator with Dirichlet boundary conditions,
$W$ is a cylindrical Wiener process in $H$, and $F$ and $R$ are defined in (\ref{heat notazioni}).
The control problem in its abstract formulation is to minimize
over all admissible controls the cost functional
\begin{equation}
J\left(  t,x,u\right)  =\mathbb{E}\int_{t}^{T} l\left(s,X^{u}_s,u_s\right)
\,ds
+\mathbb{E}\phi\left(X^{u}_T\right).
\label{heat cost abstr}%
\end{equation}
\begin{remark}\label{remarkKolmoheat}
 We notice that due to the fact that the final cost is only continuous, and the hamiltonian function
has more than linear, namely quadratic, growth with respect to $z$,
this control problem cannot be treated with techniques in the existing literature: 
in \cite{MR} lipschitz continuity of the final datum is required. In \cite{Mas1} the final datum is assumed to be bounded and continuous but the operator $G$ has to commute with the Laplace operator $A$, and this does not necessarily happens here. Hypothesis \ref{ip costo heat} ensures invertibility of $G$ which is a crucial assumption in the present paper.
\end{remark}

By applying results in section \ref{sez-appl-contr}, we get the following results.

\begin{theorem}\label{heat-th-contr}
Let $X^u$ be the solution of equation (\ref{heat equation contr}), let the cost
be defined as in (\ref{heat cost}) and let
\ref{ip costo heat} hold true.
 For every $t\in [0,T]$, $x_0\in L^2(\calo)$ and for
 all admissible control $u$ we have $J(t,x,u(\cdot))
 \geq v(t,x)$,
  and the
 equality holds if and only if
$$
u_s\in \Gamma\left( s, X^{u,t,x}_s,\nabla
v(s ,X^{u,t,x}_s)G(s)
\right)
  $$
Moreover
assume that the set-valued map $\Gamma$ is nonempty and let $\gamma$
be its measurable selection.

The closed loop equation admits a weak solution
$(\Omega,\mathcal{F},
\left(\mathcal{F}_{t}\right) _{t\geq 0}, \mathbb{P}, W,X)$
which is unique in law and setting
$$
u_{\tau}=\gamma\left(\tau, X_\tau,\nabla v(\tau,X_\tau )G(\tau)\right),
$$
we obtain an optimal admissible control system $\left(
W,u,X\right) $.
\end{theorem}

\noindent {\bf Proof.} The proof follows from the abstract formulation of the problem, and 
by applying theorems \ref{th-rel-font} and \ref{teo su controllo debole}.
\qed

\textbf{Acknowledgments}
 The author has been partially supported by the Gruppo Nazionale per l'Analisi Matematica,
la Probabilit\`a e le loro Applicazioni (GNAMPA) of the Istituto Nazionale di Alta Matematica
(INDAM).

\end{document}